\documentclass[10pt]{article}
\usepackage{amssymb,amsmath,amsfonts,mathrsfs,upgreek,textgreek,bbold} 
\usepackage{mathabx}

\usepackage[nottoc]{tocbibind}

\usepackage{graphicx} 

\usepackage[colorlinks=true]{hyperref}
\usepackage{comment}

\renewcommand\comment[1]{{\iffalse #1 \fi}}

\newtheorem{theorem}{Theorem}[section]

\newtheorem{lemma}{Lemma}[section]
\newtheorem{proposition}{Proposition}[section]

\newtheorem{assumption}{Assumption}[section]

\newtheorem{remark}{Remark}[section]
\newtheorem{notation}{Notation}[section]
\newtheorem{definition}{Definition}[section]

\newcommand{\io}{{\infty}}

\newcommand{\real}{ {\mathbb R}   }
\newcommand{\torus}{ {\mathbb T}    }
\newcommand{\integer}{ {\mathbb Z}   }
\renewcommand{\natural}{ {\mathbb N}   }
\newcommand{\complex}{ {\mathbb C}   }

\newcommand{\cB}{ {\mathcal B}   }

\newcommand{\cM}{ {\mathcal M}   }
\newcommand{\cH}{ {\mathcal H}   }

\renewcommand{\Im}{\, {\rm Im}\,}
\renewcommand{\Re}{\, {\rm Re}\,}

\newcommand{\eproof}{\qed}

\newcommand\beq[1]{ \begin{equation}\label{#1} }
\newcommand{\eeq}{ \end{equation} }
\newcommand{\beqno}{ \[ }
\newcommand{\eeqno}{ \] }
\newcommand\beqa[1]{ \begin{eqnarray} \label{#1}}
\newcommand{\eeqa}{ \end{eqnarray} }
\newcommand{\beqano}{ \begin{eqnarray*} }
\newcommand{\eeqano}{ \end{eqnarray*} }

\newcommand\dfn[1]{ \begin{definition}\label{#1} }
\newcommand\edfn{ \end{definition} }
\newcommand\ass[1]{ \begin{assumption}\label{#1} }
\newcommand\eass{ \end{assumption} }

\newcommand\notat[1]{ \begin{notation} \label{#1} 
 }
\newcommand\enotat{\end{notation}}

\newcommand\rem{\begin{remark} 
\rm 
}
\newcommand\erem{\end{remark} 
}

\newcommand{\proof}{\noindent{\bf Proof\ }}
\newcommand\equ[1]{{\rm (\ref{#1})}}
\newcommand{\nl}{{\smallskip\noindent}}
\newcommand{\giu}{{\medskip\noindent}}
\newcommand{\Giu}{{\bigskip\noindent}}

\newcommand{\qed}{\hskip.5truecm
\vrule width 1.7truemm height 3.5truemm depth 0.truemm
\par\Giu}
\newcommand{\qedeq}{\hskip.5truecm
\vrule width 1.7truemm height 3.5truemm depth 0.truemm}

\newcommand\casi[2]{ \left\{  
\begin{array}{l}
 {#1}  \\
 {#2} 
 \end{array} \right.}

\newcommand\casitwo[4]{ \left\{  \begin{array}{ll}
 {#1} & \mbox{ {\rm if} ${#2}$} \\
 {#3} & \mbox{ {\rm if} ${#4}$}
 \end{array} \right.}

\newcommand\bolla{{\tiny ${}^\bullet\,$}}

\newcommand\sopra[2]{\genfrac{}{}{0pt}{}{#1}{#2}}

\newcommand{\e}{\varepsilon}

\renewcommand{\a }{\alpha }
\renewcommand{\b }{\beta }

\newcommand{\ii }{{\rm i} }
\renewcommand{\d }{\delta }

\newcommand{\g }{\gamma}

\newcommand{\G }{\Gamma }
\renewcommand{\l }{\lambda }

\renewcommand{\t }{\tau }

\renewcommand{\O }{\Omega }
\newcommand{\C}{\mathbb{C}}
\newcommand{\Z}{\mathbb{Z}}
\newcommand{\z }{z}

\newcommand{\normadue}{|}

\newcommand{\cd}{{c_{{}_{\rm o}}}} 

\renewcommand\varpi{\upkappa}

\newcommand{\Nf}{\mathtt N}

\renewcommand\AA{{\rm A}}
\newcommand\hAA{\hat\AA}

\newcommand\ttg{{\mathtt g}}
\newcommand\ttx{{\mathtt x}}
\newcommand\tty{{\mathtt y}}
\newcommand\ttX{{\mathtt X}}
\newcommand\ttY{{\mathtt Y}}
\newcommand\ttP{{\mathtt P}}
\newcommand\ttQ{{\mathtt Q}}
\newcommand\ttq{{\mathtt q}}
\newcommand\ttp{{\mathtt p}}
\newcommand\ttG{{\mathtt G}}

\newcommand{\cgot}{\mathfrak c}

\newcommand{\Ggot}{{\mathfrak G}}

\def\N{\mathbb N}
\def\R{\mathbb R}
\def\T{\mathbb T}

\def\dst{\displaystyle}
\def\bks{\, \backslash\, }
\def\meas{{\rm\, meas\, }}
\newcommand\eqby[1]{\stackrel{\equ{#1}}{=}}

\newcommand\proiezione{\, {\mathtt p}}

\newcommand\modulo{|}
\newcommand\ham{\mathtt{H}}
\newcommand\hamk{\ham_k}

\newcommand\hamo{\ham_{\rm o}}
\newcommand\hamsec{\overline\ham}
\newcommand\hamsecu{{\overline\ham^{{}_{(1)}}}}
\newcommand\hamsecd{{\overline\ham^{{}_{(2)}}}}

\newcommand{\Hpend}{{\ham}_{\flat}}

\newcommand{\Hsharp}{{\ham}_{{}_k}}
\newcommand\Gu{{\ttG^{\!\!\scalebox{0.5}{ \rm (1)}}}}
\newcommand\Guo{{\ttG_{\rm o}^{\!\!\scalebox{0.5}{ \rm (1)}}}}
\newcommand\Gfo{\ttG_{{}_\sharp}^{\rm o}}
\newcommand\Gf{\ttG_{{}_k}}
\newcommand\bGf{\bar\ttG_{{}_k}}

\newcommand\betta{\upeta}
\newcommand\bettao{\upeta_{\rm o}}

\newcommand\lalla{\upmu}

\newcommand\sa{\theta} 
\newcommand\loge{\uplambda}

\newcommand\Gm{{\mathtt G}}
\newcommand\GO{{\bar{\mathtt G}}} 
\newcommand{\tetta}{\vartheta }
\newcommand{\tettao}{\vartheta_{\rm o}}

\newcommand\pp{{\mathfrak p}}
\newcommand\qq{{\mathfrak q}}

\newcommand\chr{\ro}

\newcommand\chs{{\check \so}}

\newcommand{\K}{{\mathtt K}}
\newcommand{\KO}{{\mathtt K}_{{}_{\rm o}}}
\newcommand{\suca}{\upepsilon}
\newcommand\morse{\upbeta}
\newcommand{\Ro}{{\mathtt R}} 
\newcommand{\ro}{{\mathtt r}}
\newcommand{\so}{{\mathtt s}}

\newcommand\act{I}  

\newcommand\Fio{\Phi_{\!{}_0}}
\newcommand\Fiuno{\Phi_{\!{}_1}} 
\newcommand\Fiunoc{\check\Phi_{\!{}_1}} 
\newcommand\Fidue{\Phi_{\!{}_2}}
\newcommand\Fiduec{\check\Phi_{\!{}_2}}
\newcommand\Fitre{\Phi_{\!{}_3}}
\newcommand\Fitrec{\check\Phi_{\!{}_3}}

\newcommand\cin{\upnu}
\newcommand\cins{\upnu_{{}_k}}

\newcommand{\ts}{\textstyle}

\newcommand{\norma}{\,\thickvert\!\!\thickvert\,}

\newcommand{\fproj}{\pi_{\!{}_{\integer k}}}

\usepackage{appendix}

\newcommand\gen{{\cal G}^n}
\def\genk{{\cal G}^n_{K}}
\def\genKO{{\cal G}^n_{\KO}}
\def\genK{{\cal G}^n_{\K}}

\renewcommand{\cgot}{\mathtt c}

\newcommand\noruno[1]{  |#1|_{{}_1} }

\renewcommand\ln{\log}
\newcommand\st{{\rm \ s.t.\ }}

\newcommand\hol{{\mathbb B}}
\newcommand{\Bns}{{\mathbb B}^n_s}
\newcommand{\Gns}{{\mathbb G}^n_{s}}
\newcommand{\Gnst}{\tilde{\mathbb G}^n_{\! s}}

\newcommand\pko{{\proiezione^\perp_k}}
\newcommand\pk{{\proiezione_k}}

\newcommand\Rz{{\mathcal R}^0}	
\newcommand\Ru{{\mathcal  R}^1}	
\newcommand\Ruk{{\mathcal  R}^{1,k}}	
\newcommand\Rd{{\mathcal  R}^2}	

\newcommand\Rzt{\widetilde{\mathcal R}^0}
	 
\newcommand\Rukt{\widetilde{\mathcal  R}^{1,k}}

\newcommand\DD{{\rm B}}
\newcommand\DDD{\mathscr D}
\newcommand\ZZ{D}

\newcommand\p{I}

\newcommand\Gdag{\Ggot_{\!{}_\dag}}

\renewcommand\diamond{{{}\!\varstar}}

\newcommand\chk{\chi_{{}_k}}

\newcommand\htk{\hat {\mathcal Q_k}}
\newcommand\hzk{\hat h_{\!{}_k}}

\newcommand\ta{{\mathtt g}}

\newcommand\giuno{{\mathtt g}_{{}_1}}  
\newcommand\gitre{{\mathtt g}_{{}_3}}

\newcommand\cdr{{c}_{{}_{\star}}}

\newcommand\itcu{c_{{{}_1}}}

\newcommand\bfco{{\bf c_{{{}_0}}}}

\newcommand\bco{{\mathtt b}_{{{}_0}}}

\newcommand\bttcd{{\mathtt c}_s}
\newcommand\itcd{c_{{{}_2}}} 

\newcommand\ttcs{{\mathtt c}_s}

\makeatletter
\newcommand{\pushright}[1]{\ifmeasuring@#1\else\omit\hfill$\displaystyle#1$\fi\ignorespaces}
\newcommand{\pushleft}[1]{\ifmeasuring@#1\else\omit$\displaystyle#1$\hfill\fi\ignorespaces}
\makeatother

\title{\bf 
Global  properties of  generic real--analytic nearly--integrable Hamiltonian systems
}

\begin{document}

\author{ 
\footnotesize L. Biasco  \& L. Chierchia
\\ \footnotesize Dipartimento di Matematica e Fisica
\\ \footnotesize Universit\`a degli Studi Roma Tre
\\ \footnotesize Largo San Leonardo Murialdo 1 - 00146 Roma, Italy
\\ {\footnotesize luca.biasco@uniroma3.it, luigi.chierchia@uniroma3.it}
\\ 
}

\maketitle


\begin{abstract}\noindent
We introduce a  new class $\Gns$ of generic real analytic potentials on $\T^n$  and study global  analytic properties of 
natural nearly--integrable Hamiltonians $\frac12 |y|^2+\e f(x)$, with potential $f\in \Gns$, on the phase space $\cM=B\times \T^n$ with $B$ a given ball in $\R^n$. The phase space $\cM$ can be covered by three sets: a `non--resonant' set, which is filled up to an exponentially small set of measure  $e^{-c \K}$ (where $\K$ is the maximal size of resonances considered)  by primary maximal KAM tori; 
a `simply resonant set' of measure $\sqrt{\e} \K^a$  and a third set of measure $\e \K^b$  which is `non perturbative', in the sense that the $\ham$--dynamics on it can be described by a natural system which is {\sl not} nearly--integrable. 
We then focus on the simply resonant set -- the dynamics of which is particularly interesting  (e.g., for Arnol'd diffusion, or the existence of secondary tori) -- and show that on such a set the secular (averaged) 1 degree--of--freedom Hamiltonians (labelled by the resonance index $k\in\integer^n$)  can be put into a universal form (which we call `Generic Standard  Form'), whose main analytic properties are controlled  by {\sl only one parameter, which is uniform in the resonance  label $k$}. 
\end{abstract}


\allowdisplaybreaks

\section*{Introduction}\label{maindefinitions}

\nl
The paper is  divided in three parts.

\nl
{\bf 1.}
In the first part we discuss generic properties of (multi--periodic) analytic functions introducing a new class $\Gns$ of functions 
 real analytic on the complex neighborhood of $\T^n$ given by 
 \beqno
 \torus^n_s:=\{x=(x_1,...,x_1)\in\complex^n:\ \ |\Im x_j|<s\}/(2\pi \integer^n)
\,.
\eeqno
 Such a class -- related but smaller than the sets $\cH_{s,\t}$ introduced in \cite{BCnonlin} -- is  generic, as
 it contains an open and dense set in the norm $\|f\|_s=\sup_k |f_k|e^{s\noruno{k}}$, and has full probability measure with respect to a  natural weighted product probability measure on Fourier coefficients.
\\
The class $\Gns$ may be  described as follows. Consider a real--analytic zero average function $f$ and consider its projection   $ \fproj f$ onto the Fourier modes proportional to a given  $k\in\integer^n\bks\{0\}$ (with components with no common divisors), which is given by 
 \beqno
\sa\in\torus\mapsto \fproj f (\sa):=\sum_{j\in\integer} f_{jk} e^{\ii j\sa}\,.
\eeqno
These one dimensional projections arise naturally, e.g.,  in averaging theory, where they are the  leading terms  of the averaged (`secular') Hamiltonians around simple resonances $\{y|\, y\cdot k= \sum y_j k_j=0\}$. 
Denoting  by  $\gen$  the set of generators of 1--dimensional maximal lattices in $\Z^n$, 
the class $\Gns$ is formed by  real--analytic zero average functions  $f$ with $\|f\|_s\le 1$, which   satisfy
\beqno
\d_{\!{}_f}=\varliminf_{\sopra{\noruno{k}\to+\io}{k\in\gen}} |f_k| e^{\noruno{k}s} \noruno{k}^n>0 \,,
\eeqno
and such that the Fourier--projection $\fproj f$ is a  Morse function with distinct critical values for all $k\in\gen$ with $\noruno{k}\le \Nf$, where $\Nf$ is a {\sl a--priori} Fourier cut--off depending only on $n,s$ and $\d_{\!{}_f}$.
\\
A remarkable feature of this class of functions is that the Fourier projection $\fproj f$ is close (in a large  analytic norm) to a  shifted rescaled cosine, 
\beqno
\fproj f(\sa)\sim |f_k|\cos (\sa+\sa_k)\,; \qquad \forall k\in\gen\,,\noruno{k}\ge \Nf\,,
\eeqno
(Proposition~\ref{pollaio} below), allowing to have uniform control of the analytic properties of  secular Hamiltonians as $\noruno{k}\to+\io$. 

\nl
 We believe that the class $\Gns$   is a good candidate  to address analytic problems in dynamical systems whenever generic results -- such as generic existence of secondary tori in nearly--integrable Hamiltonian systems\footnote{I.e.,  tori which are not a continuous deformation of integrable $(\e=0$)  flat tori; for references, see  
 \cite{Nei89}, 
 \cite{MNT}, 
 \cite{Simo}, \cite{AM},
 \cite{BCnonlin}, 
 \cite{FHL},
  \cite{BClin2}.}
 or Arnol'd diffusion\footnote{\label{AD}See 
 \cite{A64}; compare also, e.g., \cite{CG}, \cite{DH}, \cite{Z}, \cite{Ma}, \cite{T}, 
 \cite{DLS}, 
  \cite{BKZ}, 
 \cite{DS}, 
  \cite{KZ},
 \cite{GPSV},
 \cite{CdL22}, \cite{CFG}.} --  are considered. 

\nl
{\bf 2.}
In the rest of the paper, we consider natural nearly--integrable Hamiltonian systems with $n\ge 2$ degrees of freedom with Hamiltonian $\ham=\frac12 |y|^2 +\e f(x)$ ($n\ge 2$),
 with potential $f$ in the class $\Gns$ with a fixed $s>0$, on a bounded phase space $\cM=\DD\times \T^n\subset \R^n\times \T^n$; in fact, in view of the model, it is not restrictive to simply consider  $\DD=\{y\in \R^n\st  |y| <1\}$.
 \\
We, then,  introduce a covering of the action space $\DD= \Rz\cup\Ru\cup \Rd$, depending on two `Fourier cut--offs' $\K>\KO>\Nf$ ($\Nf$ as above), so that: $\Rz$ is a non--resonant set up to order $\KO$; $\Ru$ is union of neighborhoods $\Ruk$ of simple resonances 
$\{y\in \DD: y\cdot k=0\}$ of maximal order $\KO$, which are non resonant modulo $\integer k$ up to order $\K$, and $\Rd$ is a set of measure proportional to $\e \K^b$ for a suitable $b>1$ (which depends only on $n$); similar `geometry--of--resonances' analysis is typical of Nekhoroshev's theory\footnote{Compare
\cite{Nek},
\cite{BGG}, 
\cite{poschel},
\cite{DG},
\cite{Nie00},
 \cite{DH},
\cite{GCB},
 \cite{BFN}, 
 \cite{NPR}.
}.
\\
 The set $\Rd$ is {\sl a non perturbative set}, namely, it  is a set where the $\ham$--dynamics   is equivalent to the dynamics of a Hamiltonian, which is {\sl not nearly--integrable}: indeed, in the simplest non trivial case $n=2$ such a Hamiltonian is  given by $|y|^2/2 + f(x)$. 
 \\
 On the other hand the set $(\Rz\cup \Ru)\times \T^n$ is suitable for high order perturbation theory, and, 
 following the averaging theory developed in \cite{BCnonlin}, we construct high order normal forms 
 (Theorem~\ref{normalform})  so  that on $\Rz\times\T^n$ the above Hamiltonian $\ham$ is conjugated, up to an exponentially small term of $O(e^{-\KO s/3})$,  to an integrable Hamiltonian,  which depends only on action variables and it is close to $|y|^2/2$. By classical KAM theory, it then follows that this set is filled  by primary\footnote{Primary tori are smooth deformation of the flat Lagrangian integrable ($\e=0$) 
 tori.} KAM tori up to a set of measure of order  $O(e^{-\KO s/7})$. Actually, here there is a delicate point: the symplectic map realizing the above mentioned conjugation  moves the boundary of the phase space $\DD\times \T^n$ by a quantity much larger than 
 $O(e^{-\KO s/7})$, therefore, in order to get the exponentially small measure estimate on the `non--torus  set'  one needs to introduce a second covering which takes care of the dynamics close to the boundary: this is done in Lemma~\ref{sedformosa} below.
\\
The analysis on the dynamics in $\Ru\times\T^n$ is {\sl much}  more complicate. In each of the neighborhoods  $\Ruk$, which
cover the set $\Ru$  as $\noruno{k}\le \KO$, one can perform  resonant averaging theory so as to conjugate $\ham$ to still an integrable system, which however depends on the resonant angle $\ttx_1=k\cdot x$. The averaged systems with secular Hamiltonians 
$\hamsec_k(\tty,\ttx_1)$ are therefore  1D--Hamiltonian systems (one degree--of--freedom systems in the symplectic variables $(\tty_1,\ttx_1)$ depending also on  adiabatic actions 
$\tty_2,...,\tty_n$),  which are close to natural systems with   potentials $\fproj f$. Such potentials, for low $k$'s,  are rather general: for instance, they may have an arbitrary large number of separatrices depending on the particular structure of $f$. The global analytic properties of the Hamiltonians $\hamsec_k(\tty,\ttx_1)$  is the argument of the third (and main) part of this paper.

\nl
{\bf 3.}
In the third part we prove that the secular Hamiltonians  $\hamsec_k(\tty,\ttx_1)$ described in the previous item {\sl can be symplectically conjugated,  for all $\noruno{k}\le \KO$, to 1D--Hamiltonians in the standard form introduced in} \cite{BCaa23}  (see, also, Definition~\ref{morso} below).  In few words, a standard 1D--Hamiltonian (which depends on $(n-1)$ external parameters) is a one degree--of--freedom Hamiltonian system close to a natural system
with a generic potential, which may be controlled  essentially by {\sl only one parameter}, namely,  the parameter $\kappa$ appearing in Eq. \equ{alce} below; here, `essentially' means, roughly speaking,  that $\kappa$ governs  the main scaling properties of the Hamiltonian $\hamsec_k$. What is particularly relevant is that the $\kappa$ parameter of the secular Hamiltonians $\hamsec_k$  is shown to be   {\sl independent of $k$}, as it depends only on
$n$, $s$,   the above parameter $\d_{\!{}_f}$, and on a fourth parameter $\b$, which measures the Morse properties of the the potentials $\fproj f$ with $\noruno{k}\le \Nf$; compare Eq. \equ{kappa} and Remark~\ref{rampulla}--(i).\\
This uniformity allows to analyze global analytic properties: for example,   
 the action--angle map for standard Hamiltonians, as discussed  in \cite{BCaa23},  depends only on the parameter $\kappa$ of the standard Hamiltonian and therefore can be used  {\sl simultaneously} for all the secular Hamiltonians $\hamsec_k$, allowing for a nearly--integrable description of $\ham$ on $\Ruk\times\torus^n$ with uniformly exponentially small perturbations.

\nl
The results presented in this paper  may be useful  in attacking some of the fundamental open problems in the analytic theory of nearly--integrable Hamiltonian systems
such as  Arnol'd diffusion for generic real analytic systems, and provide  indispensable tools to develop a  `singular KAM Theory', namely a KAM theory dealing simultaneously with primary and secondary persistent Lagrangian tori in the full  phase space, except for the non--perturbative set $\Rd$. 
In particular,  Theorem~\ref{sivori} below is the starting point for, e.g.,  the following result, which  (up to the logarithmic correction and in the case of natural systems) proves a conjecture  by Arnold, Kozlov and Neishtadt\footnote{See, \cite[Remark~6.18, \S~6.3--C]{AKN}.},

\giu
{\bf Theorem (\cite{BClin2})}
{\sl  Fix $n\ge 2$, $s>0$,  $f\in \Gns$, $B$ an open ball in $\R^n$, and let 
$\dst \ham(y,x;\e):=\frac12 |y|^2 +\e f(x)$.
Then, there exists a constant $c>1$ such that, for all  $0<\e<1$, all  points in $B\times \T^n$
 lie on a maximal KAM torus for $\ham$, except for 
   a subset of measure bounded by
 $
 \, c\, \e |\ln \e|^{\g}$ with $\g:=11 n +4$. 
}

\giu
Let us remark that, since it is well known that the asymptotic (as $\e \to 0$) density of non--integrable {\sl primary} tori is $1-c\sqrt\e$ (see \cite{Nei}, \cite{P82}), the difference of order of the density of invariant maximal tori  in the above theorem must come from secondary tori, i.e., the tori in $\Ru\times\T^n$ whose leading dynamics is governed by the secular Hamiltonians $\hamsec_k(\tty,\ttx_1)$ discussed in this paper. 

\section{Generic real analytic periodic functions} \label{gennaro}

We begin with a few definitions.

\dfn{pelle}{\bf (Norms on real analytic periodic functions)}\\
For $s>0$ and  $n\in\natural=\{1,2,3...\}$,
   consider  the Banach space of
zero average
real analytic periodic functions  $\dst f:x\in \torus^n:=\real^n/(2\pi \Z^n)\mapsto \sum_{  k\in\integer}f_k e^{\ii k\cdot x}$, $f_0=0$,  
with finite norm\footnote{As usual $\noruno{k}:=\sum |k_j|$.}       
\beqno
 \|f\|_s:=\sup_{k\in \integer^n}  |f_k| e^{\noruno{k}s}\,,
\eeqno
and denote by $\Bns$ the closed unit ball of functions $f$ with $\|f\|_s\le 1$.
\edfn

\nl
Besides the norm $\|\cdot\|_s$, we shall also use the following  two  (non equivalent) norms
$$
{\modulo}f{\modulo}_{s}:=\sup_{\T^n_s}|f|\,,\qquad 
\norma f\norma_{s}:=
\sum_{k\in\Z^n} 
|f_k| e^{|k|_{{}_1}s}\,.
$$
Such norms satisfy the  relations
 $$\|f\|_{s}\leq {\modulo}f{\modulo}_{s}
\leq
\norma f\norma_{s}\,.
$$
Notice also the following 
 `smoothing property'  of the norm  $\norma \cdot\norma_{r,s}$:
{\sl if $s'\leq s$, then for any $N\ge 1$, one has} 
\begin{equation}\label{lesso}
f(y,x)=\sum_{\noruno{k}\geq N}f_k(y)e^{\ii k\cdot x}
\qquad
\Longrightarrow
\qquad
\norma f \norma_{r,s'}\leq e^{-(s-s')N} \norma f \norma_{r,s}\,.
\end{equation}

\dfn{generators}{\bf (Generators and Fourier projectors)}\\
{\rm (i)} 
Let $ \integer^n_\varstar$ be the set of integer vectors $k\neq 0$ in $\integer^n$ such that the  first non--null  component is positive:
\beq{iguazu}
 \integer^n_\varstar:=
 \big\{ k\in\integer^n:\ k\neq 0\ {\rm and} \ k_j>0\ {\rm where}\ j=\min\{i: k_i\neq 0\}\big\}\,,
 \eeq
and denote by $\gen$ 
the set of {\sl generators of 1d maximal lattices} in $\integer^n$, namely, the set of  vectors $k\in  \integer^n_\varstar$ 
such that the greater common divisor ({\rm gcd})  of their components is 1:
\beqno
\gen:=\{k\in \integer^n_\varstar:\ {\rm gcd} (k_1,\ldots,k_n)=1\}\,.
\eeqno
Let us also denote  by $\genk$ the generators of size not exceeding $K\ge 1$,
\beqno
\genk:=\gen \cap \{\noruno{k} \leq K \}\,,
\eeqno
{\rm (ii)} 
Given  a zero average  real analytic 
periodic function 
and $k\in \gen$,
we define
 \beq{kproj}
\sa\in\torus\mapsto \fproj f (\sa):=\sum_{j\in\integer} f_{jk} e^{\ii j\sa}\,.
\eeq
\edfn
Notice that $f$ can be uniquely written as 
\beqno
f(x)=
\sum_{k\in \gen} \fproj f (k\cdot x)\,.
\eeqno

\begin{definition}\label{buda}
Let $\b>0$. A function $F\in C^2(\T,\R)$   is called a   {\bf $\b$--Morse function}   if 
\beqno
\min_{\sa\in\T} \ \big( |F'(\sa)|+|F''(\sa)|\big) 
\geq\b \,,
\quad
\min_{i\neq j } 
|F(\sa_i)-F(\sa_j)|
\geq \b \,, 
\eeqno
where $\sa_i\in\T$ are the  critical points of $F$.
\end{definition}

\dfn{pigro}{\bf (Cosine--like functions)}
 Let  $0<\ttg< 1/4$.
 We say that a real analytic function  $G:\T_1\to\C$  is
$\ttg$--cosine--like
if, for  some $\eta>0$ and 
 $\sa_0\in\R$, one~has 
  \beqno
 {\modulo}G(\sa)-\eta\cos (\sa+\sa_0) {\modulo}_1
 := \sup_{|\Im \sa|<1}{\modulo}G(\sa)-\eta\cos (\sa+\sa_0) {\modulo}
 \leq \eta\ttg\,.
\eeqno
 \edfn
Notice that this notion is invariant by rescalings: $G$ is $\ttg$--cosine--like if and only if $\l G$ is $\ttg$--cosine--like for any $\l>0$. Beware of the usage of $|\cdot|_1$ as sup norm on  $\torus_1$, the complex strip of width~2.

\nl
Now, the main definition. 

 \dfn{sicuro} {\bf (The analytic class $\Gns$)}
We  denote by $\Gns$ the subset  of functions $f\in \Bns$  such that the following two properties 
hold:
\beqa{P1}
&&
\varliminf_{\sopra{\noruno{k}\to+\io}{k\in\gen}} |f_k| e^{\noruno{k}s} \noruno{k}^n>0\,,
\\
&& 
 \forall \ k\in\gen\,,\ \fproj f\ {\rm is \ a\ Morse\ function\ with \ distinct \ critical\ values}\,.
 \nonumber
\eeqa
\edfn

\rem\label{posizione}
(i)
If $f\in  \Bns$, then the function $\fproj f$ belongs to $ \hol_{|k|_{{}_1}s}^1$ and therefore has a domain of analyticity which increases with the norm of $k$.

\nl
(ii) A simple example of function    in $\Gns$ is given by
\beqno
f(x):=2 \sum_{k\in\gen} e^{-s\noruno{k}} \cos k\cdot x\,.
\eeqno
Indeed, one checks immediately that
\beqno
\|f\|_s=1\,,\qquad 
\varliminf_{\sopra{\noruno{k}\to+\io}{k\in\gen}} |f_k| e^{\noruno{k}s} \noruno{k}^n=+\io\,,\qquad
\fproj f (\sa)=2 e^{-s\noruno{k}} \cos \theta\,.
\eeqno
(iii) 
The  critical points of an analytic  Morse function on $\torus$, by compactness, cannot accumulate, hence, there  are a finite, even  number of them, which are, alternately, a  relative strict maximum and a relative strict minimum. 
In particular, if $G$ is $\b$--Morse, then the number  of its critical points  can be bounded by
$\pi\sqrt{2\max |G''|/\b}$. 
Indeed, if $\sa\neq \sa'$  are critical points of $G$, then, by \equ{cimabue}  one has
$$\ts\b\le|G(\sa)-G(\sa')|\le \frac12 (\max|G''|) \,{\rm dist}(\sa,\sa')^2\,,$$ which implies that the minimal distance between two critical points is  $\sqrt{2\b/\max|G''|}$ and the claim follows.
\erem

\subsection{Uniform behaviour of large-mode Fourier projections}\label{ironman}

If a function $f\in \Bns$ satisfies \equ{P1}, then, {\sl apart from a finite number of Fourier modes,  its Fourier projections
$\fproj f$ are close to a shifted rescaled cosine}, a fact that allows, e.g.,  to have a uniform analytic theory of high order perturbation theory.

\nl
To discuss this matter, let us first point out that
for any sequence of real numbers $\{a_k\}$ and for any  function $N(\d)$ such that  
$\lim_{\d\downarrow 0} N(\d)=+\infty$ one has 
\beq{analisi1}
\varliminf a_k>0 \quad \iff \quad \exists \ \d>0\ \ {\rm s.t.}\ \ a_k\ge \d\,,\ \forall\ k\ge N(\d)\,.
\eeq 
We shall apply this  remark to the minimum limit in \equ{P1} with a particular choice of the function $N(\d)$, namely,  we define
$\Nf(\d)=\Nf(\d;n,s)$ as 
\beq{enne}
\Nf(\d):=2\, \max\Big\{1\,,\,\frac1s \,   \log \frac\cd{s^n\, \d}\Big\} \,,\qquad \cd:= 2^{44}\ (2n/e)^n\,.
\eeq 
For later use, we point out that\footnote{In fact, if $s\ge 1$ then $\Nf\ge 2\ge 2/s$, while if $s<1$ then the logarithm in
\equ{enne} is larger than one, so that $\Nf\ge 2/s$ also in this case.}  
\beq{bollettino1}
\Nf\ge2  \ttcs\,,\quad {\rm where}\quad \ttcs:=\max\big\{1, 1/s\big\}\,.
\eeq
From \equ{analisi1} it follows  that if $f$ satisfies \equ{P1},  one can find  $0<\d\le 1$ such that
\beq{P1+}
|f_k|\geq \d \noruno{k}^{-n}\, e^{-\noruno{k} s}\,,\qquad \forall \ k\in\gen\,,\ \noruno{k}\ge\Nf\,.
\eeq 
The main feature of the above choice of $\Nf$ is    that, for $\noruno{k}\ge\Nf$, $\fproj f$ is very close to a shifted rescaled cosine function:

\begin{proposition}\label{pollaio}  
Let $\d>0$, $f\in\Bns$ and  assume \equ{P1+}.
 Then, 
 for any  $k\in \gen $ with  $ \noruno{k}\geq \Nf$,   
 $\fproj f$ is $2^{-40}$--cosine--like (Definition~\ref{pigro}).  
\end{proposition}
\proof 
We shall prove something slightly stronger, namely,  that there exists $\sa_k\in[0,2\pi)$ 
so that 
\begin{equation}\label{alfacentauri}
\fproj f(\sa)=2 |f_k| \big(\cos(\sa+ \sa_k)+F^k_\star(\sa)\big)\,,\quad F^k_{\! \varstar}(\sa):=\frac{1}{2|f_k|}\sum_{|j|\geq 2}f_{jk}e^{\ii j \sa}\,,
\end{equation}
with $F^k_{\! \varstar}\in\hol_1^1$ and (recall the definition of the norms in \equ{norme})
 \begin{equation}\label{gallina}
\modulo F^k_{\! \varstar} \modulo_1\le \norma F^k_{\! \varstar}\norma_{1}\leq 2^{-40}\,. 
\end{equation}
Indeed,  by definition of $\fproj f$, 
\beqno
\fproj f(\sa):= \sum_{j\in\integer\bks\{0\}} f_{jk} e^{\ii j \sa}
= \sum_{|j|=1} f_{jk}e^{\ii j\sa} + \sum_{|j| \ge 2} f_{jk}e^{\ii j\sa} \,,
\eeqno
and, defining $\sa_k\in[0,2\pi)$ so that $e^{\ii \sa_k}= f_k/|f_k|$, one has
$$
\frac{1}{2|f_k|}\sum_{|j|=1}f_{jk}e^{\ii j \sa}=\Re \Big( \frac{f_k}{|f_k|} e^{\ii \sa}\Big)=\Re e^{\ii (\sa+\sa_k)}
=\cos (\sa+\sa_k)\,,
$$
which yields   \equ{alfacentauri}. Now, since $f\in\Bns$ it is 
$|f_k|\le e^{-\noruno{k}s}$
 and, by  \equ{P1+},   $|f_k|\geq \d \noruno{k}^{-n}\, e^{-\noruno{k} s}$.
 Therefore, for $\noruno{k}\ge \Nf$, one has 
\beqa{onlyyou}
\norma F^k_{\! \varstar}\norma_{{}_1}&\stackrel{\equ{alfacentauri}}=&
\frac{1}{2|f_k|}\sum_{|j|\geq 2}|f_{jk}|e^{|j|}
\leq
\frac{\noruno{k}^n e^{\noruno{k}s}}{2\d}\sum_{|j|\geq 2}|f_{jk}|e^{|j|}\nonumber
\\
&\le&
\frac{\noruno{k}^n e^{\noruno{k}s}}{2\d}\sum_{|j|\geq 2}e^{-|j|(\noruno{k}s-1)}
\nonumber
\\
&\le&
\frac{2 e^2 \noruno{k}^n}{\d} \ e^{-\noruno{k}s}
=\frac{2^{n+1} e^2 }{s^n\d} e^{-\frac{\noruno{k}s}2}\ \ \Big(\frac{\noruno{k}s }2\Big)^n e^{-\frac{\noruno{k}s}2}
\nonumber
\\
&\le& \Big(\frac{2n}{e s}\Big)^n\, \frac{2e^2}{\d} \,  e^{-\frac{\Nf s}2}
\le 2^{-40}\,,
\eeqa
where the geometric series converges since $\noruno{k}s\ge \Nf s\ge2 $ (by \equ{bollettino1}) and   last inequality follows 
by  definition of $\Nf$ in \equ{enne}.
\qed

\rem
In fact, the particular form of $\Nf$ is used {\sl only} in the last inequality in \equ{onlyyou}. 
\erem

\noindent
Next, we need an elementary calculus lemma:

\begin{lemma}\label{pennarello}
 Assume that $F\in C^2(\T,\R)$, $\bar\sa$  and $0<\cgot<1/2$ are such that
  $$\|F-\cos (\sa+\bar \sa)\|_{C^2}\le \cgot\,,$$
  where  $\|F\|_{C^2}:=\max_{0\leq k\leq 2}\sup|F^{(k)}|$. 
 Then,  $F$ has only two critical points and it is $(1-2 \cgot)$--Morse (Definition~\ref{buda}).
 \end{lemma}
\proof By considering the translated function $\sa\to F(\sa-\bar\sa)$, one can reduce oneself to the case $\bar\sa=0$ ($F$ is $\b$--Morse, if and only if $\sa\to F(\sa-\bar\sa)$ is $\b$--Morse).\\
Thus, we set $\bar \sa=0$, and note that, by assumption $|F'|=|F'+\sin\sa-\sin\sa|\ge|\sin \sa|-\cgot $, and, analogously, $|F''|\ge |\cos\sa|-\cgot $. Hence, $|F'|+|F''|\ge|\sin\sa|+|\cos\sa|-2\cgot \ge 1-2\cgot $. Next, let us show that $F$ has a unique strict maximum $\sa_0\in I:=(-\pi/6,\pi/6)$ (mod $2\pi$). Writing $F=\cos\sa+g$, with $g:=F-\cos \sa$, one has that $F'(-\pi/6)=1/2+g'(\pi/6)\ge 1/2 - \cgot >0$, and, similarly $F'(\pi/6)\le -1/2 +\cgot $, thus $F$ has a critical point in $I$, and, since $-F''=\cos\sa -g''\ge \cos\sa-\cgot \ge \sqrt3/2-\cgot >0$, $F$ is strictly concave in $I$, showing that such critical point is unique and it is a strict local minimum. In fact, similarly one shows that $F$ has a second critical point $\sa_1\in (\pi-\pi/6,\pi+\pi/6)$ where $F$ is strictly convex, so that $\sa_1$ is a strict local minimum; but, since
 in the complementary of these intervals $F$ is strictly monotone (as it is easy to check),  it follows that $F$ has a unique global strict maximum and a unique global strict minimum. 
Finally, $F(\sa_0)-F(\sa_1)\ge \sqrt3-2\cgot >1-2\cgot $ and the claim follows. \qed

\nl
From  Proposition~\ref{pollaio} and Lemma~\ref{pennarello} one gets immediately:

\begin{proposition}\label{punti} Let $\d>0$, $f\in\Bns$ and  assume \equ{P1+}.
Then,
for every
$k\in \gen$ with  $\noruno{k}\ge \Nf$, the function
 $\fproj f$ is $|f_k|$--Morse.
\end{proposition}
\proof   As in the proof of Proposition~\ref{pollaio},
we get
\beq{derby}
\Big|\frac{\fproj f}{2f_k} -  \cos(\sa+\sa^k)\Big|_1\stackrel{\equ{alfacentauri}}=
 |F^k_\star|_1\leq
 \norma F^k_\star\norma_1\stackrel{\equ{gallina}}
\leq 2^{-40}\,,
\eeq
which implies that the function $F:=\fproj f/(2f_k)$ is $C^2$--close to  a (shifted) cosine: Indeed, by Cauchy estimates
$\|\cdot\|_{C^2}\leq 2 |\cdot|_1$, so that 
\beqno
\|F-\cos(\sa+\sa^k)\|_{C^2}=\max_{0\le j\le 2} \max_\T |\partial_\sa^j(F-\cos(\sa+\sa^k))|\le 
2|F^k_\star|_1 \stackrel{\equ{derby}} \leq  2^{-39} \,.
\eeqno
By Lemma~\ref{pennarello} we see that $F$ is $(1-2^{-38})$--Morse, and the claim follows by rescaling.~\qed

\subsection{Genericity}

In this section we prove that $\Gns$ is a generic set in  $\Bns$.

\dfn{sicuro2} Given $n,s>0$, $0<\d\le 1$ and $\b>0$ and  $\Nf$ as in \equ{enne} we call 
$\Gns(\d,\b)$ the set of functions in $\Bns$ which satisfy \equ{P1+} together with:
 \beq{P2+}
\fproj f\ {\rm is \ \b\!-\!Morse}\,,\qquad \, \ \ \forall \ k\in\gen\,,\ \noruno{k}\le \Nf\,.
\eeq
\edfn
Then, the following lemma holds:

\begin{lemma}
 \label{telaviv} Let $n,s>0$. Then, 
$\dst
\Gns=\bigcup_{\sopra{\d\in (0,1]}{\b>0}}\Gns(\d,\b)$.
\end{lemma}

\proof
Assume $f\in \Gns$  and 
let $0<\d_0\le 1$ be smaller than limit inferior in \equ{P1}.
Then, there exists $N_0$ such that $|f_k|>\d_0 \noruno{k}^{-n} e^{-\noruno{k}s}$, for any $\noruno{k}\ge N_0$, $k\in\gen$.
Since $\lim_{\d\to 0} \Nf=+\infty$, there exists $0<\d<\d_0$ such that $\Nf>N_0$. Hence, if $\noruno{k}\ge \Nf$ and 
$k\in\gen$, 
\equ{P1+} holds.
\\
Since $\fproj f$ is, for any $\noruno{k}\le \Nf$, a Morse function with distinct critical values one can, obviously, find a $\b>0$ for which
\equ{P2+} holds. Hence $f\in \Gns(\d,\b)$.

\nl
Now, let $f\in \bigcup \Gns(\d,\b)$. Then, there exist $\d\in (0,1]$ and $\b>0$ such that \equ{P1+} and \equ{P2+} hold. Then, \equ{P1}   follows immediately from \equ{P1+}. 
By Proposition~\ref{pollaio}, for any $k\in \gen$ with $\noruno{k}> \Nf$, 
  $\fproj f$ is $2^{-40}$--cosine--like, showing (Lemma~\ref{pennarello}) that $\fproj f$ is Morse with distinct critical values also for $\noruno{k}\ge \Nf$. The proof is complete.
\qed

\begin{proposition}\label{adso}
$\Gns$ contains an open and dense set in $\Bns$. 
\end{proposition}

\nl
To prove this result we need a preliminary  elementary result on real analytic periodic functions:

\begin{lemma}\label{trifolato} Let $F=\sum F_j e^{\ii j \sa}$ be a real analytic function on $\torus$. 
There exists a compact set $\G\subseteq\C$ (depending on $F_j$ for $|j|\ge 2$) of zero Lebesgue  measure such that if the Fourier coefficient $F_1$  does not belong to $\G$, then 
$F$ is a Morse function with distinct critical values.
\end{lemma}
\proof
Without loss of generality we may assume that $F$ has zero average. Then, letting $z:=F_1\in\complex$, we write $F$ as
\begin{equation}\label{efesta}
F(\sa)=
z e^{\ii \sa} + \bar z e^{-\ii \sa} + G(\sa):=z e^{\ii \sa} + \bar z e^{-\ii \sa} +
\sum_{|j|\geq 2} F_je^{\ii j\sa} \,.
\end{equation}
 When $G\equiv 0$ the claim is true with   $\G=\{0\}$.\\
Assume that  $G\not\equiv 0$. 
Observe that, since $G$ is real--analytic, the equations 
$ F'(\sa)=0= F''(\sa)$ are  equivalent to the single equation
$z=\frac12 e^{-\ii \sa} \big( \ii G'(\sa)+G''(\sa) \big)$, which, as $\sa\in\torus$, describes a smooth closed `critical' curve $\Gamma_1$ in $\complex$.
\\
Observe also that $F$ has distinct critical points $\sa_1,\sa_2\in\T$ with the same critical values if and only if
the following three real equations are satisfied:
\beq{odessa}
F'(\sa_1)=0\,,\qquad
F'(\sa_2)=0\,,\qquad
F(\sa_1)-F(\sa_2)=0\,.
\eeq
We claim that if $z,\sa_1,\sa_2$ satisfy \equ{odessa}  
then
 \begin{equation}\label{celebration}
  z=\zeta(\sa_1,\sa_2)\,,\qquad
 g(\sa_1,\sa_2)=0\,,
\end{equation}
 with $\zeta$ and $g$ real analytic on $\T^2$ given by
 \beqano
&&\zeta(\sa_1,\sa_2):=
\left\{\begin{array}l
\frac{\ii}{2(e^{\ii \sa_1}-e^{\ii \sa_2})}
\big(
G'(\sa_1)-G'(\sa_2) +\ii G(\sa_1) -\ii G(\sa_2)
\big)\,,\quad \mbox{for}\ \  \sa_1\neq\sa_2\,;
\\ \ \\
\frac1{2e^{\ii \sa_1}}
\big(
G''(\sa_1)+\ii G'(\sa_1)
\big)\,,\phantom{AAAAAAAAAAAAAaa} \mbox{for}\ \   \sa_1=\sa_2\,,
\end{array}\right.
\\ \ \\
&&
 g(\sa_1,\sa_2)
:=
\big(1-\cos(\sa_1-\sa_2)\big)\big( G'(\sa_1)+G'(\sa_2) \big)
- \sin (\sa_1-\sa_2)\big( G(\sa_1)-G(\sa_2) \big)\,.
\eeqano
Indeed, summing up the the third equation in \eqref{odessa} with 
the difference of the first two equations 
multiplied by $-\ii$, 
we get
$$
2(e^{\ii \sa_1}-e^{\ii \sa_2})
z-\ii
\big(
G'(\sa_1)-G'(\sa_2) +\ii G(\sa_1) -\ii G(\sa_2)
\big)=0\,,
$$
which is equivalent to
$z=\zeta(\sa_1,\sa_2)$. Then, by definition $g(\sa_1,\sa_1)=0$, while if $\sa_1\neq \sa_2$, 
substituting $z=\zeta(\sa_1,\sa_2)$
in the first  equation in \eqref{odessa}
and multiplying by  
$1-\cos(\sa_1-\sa_2)$
we get $g(\sa_1,\sa_2)=0$ also for $\sa_1\neq\sa_2$. Thus,  \eqref{celebration} holds.
\\
Next, we claim that the real analytic function 
$ g(\sa_1,\sa_2)$ 
is not identically zero.  
Assume by contradiction that $g$ is identically zero.
Then $g(\sa_2+t,\sa_2)\equiv 0$ for every $\sa_2$ and 
$t$, and taking the fourth derivative with respect to $t$ evaluated at  
$t=0$, we see that 
$
G'''(\sa_2)+G'(\sa_2)=0$, for all $\sa_2$.
The general (real) solution of the such differential equation 
is  given by $G(\sa_2)= c e^{\ii \sa_2} + \bar c e^{-\ii \sa_2}+c_0,$
with $c\in\C,$ $c_0\in\R, $
which contradicts  the fact that, by definition, $G_j=0$ for $|j|\le 1$. 
Thus, $ g(\sa_1,\sa_2)$ 
is not identically zero
and, therefore, 
the set
$\mathcal Z\subseteq\T^2$ of its zeros  is compact and has zero Lebesgue  measure\footnote{Compare, e.g., Corollary 10, p. 9 of \cite{GR}.}.
Clearly, also the set $\Gamma_2:=\zeta(\mathcal Z)\subseteq\C$
is compact and has zero measure, and, therefore, if we define
$\G=\Gamma_1\cup\Gamma_2$, we see that the lemma holds also in the case $G\nequiv 0$.
\qed

\proof {\bf of Proposition~\ref{adso}}
Let  $\Gnst(\d,\b)$ denote the subset of functions in $\Gns(\d,\b)$ satisfying the (stronger) condition\footnote{Here, we explicitly indicate  
the dependence on $\d$, while $n$ and $s$ are fixed.
Recall that $\Nf(\d)$ is decreasing.}
\beq{starstar}
|f_k|> \d \, e^{-\noruno{k} s}\,,\qquad \forall \ k\in\gen\,,\ \noruno{k}\ge\Nf
=\Nf(\d)\,,
\eeq
and let
$\dst
\Gnst=\bigcup_{\stackrel{0<\d\le 1}{\b>0}} \Gnst(\d,\b)\,.
$
We claim that $\Gnst$ is an open subset of  $\Bns$.
 Let $f\in\Gnst(\d,\b)$ for some 
${0<\d\le 1},{\b>0}$ and let us show 
that there exists $0<\d'\leq\d/2$ such that  if 
$g\in \Bns$ with
$\|g-f\|_s<\d'\leq \d/2$, 
then 
$g\in\Gnst(\d',\b')$ with $\b':=\min\{\b,\d e^{-s\Nf(\d/2)}\}/2$.
Indeed
$$
|\tilde f_k|e^{|k|_{{}_1}s}\geq |f_k| e^{|k|_{{}_1}s} -\|g-f\|_s > \d-\d'\geq\d/2\,,
\qquad \forall \ k\in\gen\,,\ \noruno{k}\ge\Nf(\d)\,, 
$$
namely $g$ satisfies  \equ{starstar} with $\d/2$ instead of $\d$.
We already know that
$\fproj f$ is $\b\!-\!$Morse $ \forall \ k\in\gen,\, \noruno{k}<\Nf(\d)$.
Moreover, by Proposition~\ref{punti},
we know that $\fproj f$ is $|f_k|$--Morse for
$k\in \gen$ with  $\noruno{k}\ge \Nf(\d)$.
In conclusion, by \eqref{starstar}, we get that 
$\fproj f$ is $2\b'\!-\!$Morse $ \forall \ k\in\gen,\, \noruno{k}<\Nf(\d/2)$.
Since the $\|\cdot\|_s$--norm  is stronger than the $C^2$--one, taking $\d'$ small enough we get that  
$\fproj g$ is $\b'\!-\!$Morse $ \forall \ k\in\gen,\, \noruno{k}<\Nf(\d/2)$.

\nl
Let us now show that $\Gnst$ is dense  in $\Bns$.
Fix $f$ in $\Bns$ and $0<\loge <1$. We have to find $g\in\Gnst$
with $\|g-f\|_s\leq \loge $.
Let $\d:=\loge /4$ and denote by  $f_k$ and $g_k$ (to be defined) be the Fourier coefficients of, respectively,  
$f$ and $g$. 
It is enough to define $g_k$ only for 
$k\in\integer^n_\varstar$ since, for $k\in -\integer^n_\varstar$
we set $g_k:=\bar{g}_{-k}$, since $g$ must be real analytic. 
Set $g_k:=f_k$ for
 $k\in\integer^n_\varstar\setminus\gen$.
 For
 $k\in\gen$, $|k|_{{}_1}\geq \Nf(\d)$,
 we set  $g_k:=f_k$ if
  $|f_k|e^{|k|_{{}_1}s}>\d$ and
 $g_k:= 2\delta e^{-|k|_{{}_1}s}$ otherwise.
 Consider now $k\in\gen$, $|k|_{{}_1}< \Nf(\d)$.
 We make use of Lemma~\ref{trifolato}
 with $F=\fproj g$, $z=F_1=g_k$.
Thus, by Lemma~\ref{trifolato},  there exists  a compact set $\G_k\subseteq\C$ (depending on $F_k$ for $|k|\ge 2$) of zero measure such that if $g_k\notin \G_k$ the function 
$\fproj g$ is a Morse function with distinct critical values.
We conclude the proof of the density choosing
$|g_k|<e^{-|k|_{{}_1}s}$, $|f_k-g_k|\leq \loge
e^{-|k|_{{}_1}s}$ with $g_k\notin \Gamma_k$.
\qed

\subsection{Full   measure}
Here we show that {\sl $\Gns$ is a set of probability 1 with respect to the standard product probability measure on $\Bns$}.
More precisely,  consider  the  space\footnote{$\integer^n_\varstar$ was 
defined in \eqref{iguazu}.} 
${\rm D}^{{\integer^n_\varstar}}$, where ${\rm D}:=\{w\in\complex:\ |w|\le 1\}$, endowed with the  product topology\footnote{
By  Tychonoff's Theorem,
${\rm D}^{{\integer^n_\varstar}}$ with the product topology is a
compact Hausdorff space. 
}.
The product $\sigma$-algebra of the Borel sets of ${\rm D}^{{\integer^n_\varstar}}$ is the $\sigma$--algebra generated by 
the cylinders $\bigotimes_{k\in{\integer^n_\varstar}} A_k$, where $A_k$ are
Borel sets of ${\rm D}$, which differs from ${\rm D}$ only for a finite number of $k$.
The  probability product measure $\mu_{{}_\otimes}$ on 
 ${\rm D}^{{\integer^n_\varstar}}$
is then defined by letting
$$
\mu_{{}_\otimes} \big(\bigotimes_{k\in{\integer^n_\varstar}} A_k \big):=
\prod_{k\in{\integer^n_\varstar}} 
 |A_k |\,,
 $$
where $|\cdot|$ denotes  the normalized ($|{\rm D}|=1$) Lebesgue measure on ${\rm D}$.
The (weighted) Fourier bijection\footnote{$f$ is  real analytic so that $f_{-k}=\bar f_k.$}
\begin{equation}\label{odisseo}
\mathcal F:f(x)=\sum_{k\in \integer^n_\varstar}
f_k e^{\ii k\cdot x}+\bar f_k e^{-\ii k\cdot x}
\in  \Bns \to \big\{f_k e^{|k|_{{}_1}s}\big\}_{k\in {\integer^n_\varstar}}\in \ell^\io({\integer^n_\varstar})
\end{equation}
induces a product topology on $\Bns$
and a {\sl  probability product measure} $\mu$ on the
product $\sigma$-algebra
$\cB$ of the Borellians  in $\Bns=\mathcal F^{-1}
\big({\rm D}^{{\integer^n_\varstar}}\big)$
(with respect to the induced product topology), i.e.,  given $B\in\cB$, we set $\mu(B):=\mu_{{}_\otimes}(\mathcal F(B))$.  Then one has:

\begin{proposition}\label{melk}
 $\Gns\in \cB$ and $\mu(\Gns)=1$.
 \end{proposition}

\proof  
 First note that, for every $\d,\b>0$ the set
 $\Gns(\d,\b)$ is
 closed  with respect to the product topology.
 Indeed for every $ k\in\gen$
 the set $\{f\in \Bns\ \mbox{s.t.}\ 
 |f_k|\geq \d \noruno{k}^{-n}\, e^{-\noruno{k} s}\}$
 is a closed cylinder. Moreover also 
 the set $\{f\in \Bns\ \mbox{s.t.}\ \fproj f \ \mbox{is}\ \b\mbox{--Morse}\}$
 is closed w.r.t the product topology. In fact we prove that
 the complementary $E:=\{f\in \Bns\ \mbox{s.t.}\ \fproj f \ \mbox{is not}\ \b\mbox{--Morse}\}$ is open  w.r.t the product topology.
 Indeed if $f^*\in E$  there exists a $r>0$ small enough such that
 $E_r:=\{f\in \Bns\ \mbox{s.t.}\ \|\fproj f-\fproj f^*\|_{C^2}<r \}\subseteq E.$
 Define the open cylinder
 $$
 E_{\rho,J}:=\{
 f\in \Bns\ \mbox{s.t.}\ 
|f_{jk}-f_{jk}^*|<\frac{\rho}{\noruno{j}^2}e^{-\noruno{jk}s}\ \mbox{for}\ 
j\in\Z\,,\ 
0<\noruno{j}\leq J \}\,.
 $$
 We claim that $E_{\rho,J}\subseteq E_r$ for suitably small $\rho$
 and large $J$ (depending on $r$ and $s$).
 Indeed, when $f\in E_{\rho,J}$
 $$
 \|\fproj f-\fproj f^*\|_{C^2}
 \leq
 3\sum_{j\neq 0} \noruno{j}^2 |f_{jk}-f_{jk}^*|
 \leq 
 3\rho\sum_{0<\noruno{j}\leq J}e^{-\noruno{jk}s}
 +
 6\sum_{\noruno{j}> J}\noruno{j}^2 e^{-\noruno{jk}s}
 <r
 $$
 for suitably small $\rho$
 and large $J$.
 Therefore $E_{\rho,J}\subseteq E_r\subseteq E$ and $E$ is open in
 the product topology.
 In conclusion, taking the intersection over $k\in\gen$,
 we get that  $\Gns(\d,\b)$ is
 closed  with respect to the product topology.
\\
Recalling Lemma~\ref{telaviv}, we note that $\Gns$ can be written as 
$\dst
\Gns=\bigcup_{h\in\N} \Gns(1/h,1/h)$. 
Thus $\Gns$ is Borellian.

\giu
Let us now prove that $\mu(\Gns)=1$.
Fix  $0<\d\le 1$ and denote by  
$\Gns(\d)$ be the subset of functions in $\Bns$ satisfying \equ{P1+}
and such that $\fproj f$
 is  a Morse function with distinct  critical values
for every $ k\in\gen$.
Recall \eqref{odisseo} and define
$$
\mathbb P_\d:=\mathcal F(\Gns(\d))\subseteq \ell^\io({\integer^n_\varstar})\,.
$$
Fix $\hat g=(g_k)_{k\in\integer^n_\varstar\setminus\gen}\in
\ell^\io({\integer^n_\varstar\setminus\gen})$ with $|g_k|\leq 1$ for every
$k\in\integer^n_\varstar\setminus\gen$.
  Consider
the section
\beqno
\mathbb P_\d^{\hat g}:=\big\{ \check g=(g_k)_{k\in\gen},\ |g_k|\leq 1 \ \ \mbox{s.t}\ \ |g_k|\geq \d \noruno{k}^{-n} \ \mbox{if}\ \noruno{k}\geq \Nf\,,\ \ 
g_k e^{-\noruno{k}s}\notin \G_k\,, \ \mbox{if}\ \noruno{k}< \Nf
\big\},
\eeqno
where the sets $\G_k$ (depending on $\hat g$)
   were defined in the proof of Proposition \ref{adso}
  so that, for every $k\in\gen,$ $|k|_{{}_1}< \Nf$, 
  if $g_k e^{-\noruno{k}s}\notin \G_k$ then the function\footnote{Recall \eqref{efesta}.}
$$
 g_k e^{-\noruno{k}s} e^{\ii \sa} + \bar g_k e^{-\noruno{k}s} e^{-\ii \sa} + 
 \sum_{|j|\geq 2} \hat g_{jk} e^{-\noruno{jk}s} e^{\ii j\sa}
 =\fproj f\,,\ \ \mbox{with}\ \ f:=\mathcal F^{-1}(g)\,,\ \ g=(\check g,\hat g)\,,
$$
 is a Morse function with distinct critical values.   Then, since every $\G_k$ has
 zero measure 
$$
\mu_{{}_\otimes}|_{\ell^\io({\gen})}(\mathbb P_\d^{\hat g})=
\prod_{k\in\gen, |k|_{{}_1}\geq \Nf} (1-
  \d^2\, |k|_{{}_1}^{-2{n}})\geq 1-c\d^2\,,
$$
for a suitable constant $c=c(n)$.
Since the above estimate holds for every
$\hat g\in
\ell^\io({\integer^n_\varstar\setminus\gen})$,
 by Fubini's Theorem we get
$$
 \mu_{{}_\otimes}|_{\ell^\io({\gen})}(\mathbb P_\d^{\hat g})=
  \mu_{{}_\otimes}(\mathbb P_\d)=\mu(\Gns(\d))
\geq 1-c\d^2\,.
$$
Then,
$$
\mu(\Gns)=\lim_{\d\to 0^+} \mu(\Gns(\d))=1\,. \qedeq
$$

\section{Averaging, coverings and normal forms}

In the rest of the paper we consider
{\sl real--analytic, nearly--integrable natural Hamiltonian systems}
\beqa{ham}\ts
&&\left\{
\begin{array}{l}
\dot y = -\ham_x(y,x)\\
\dot x= \ham_y(y,x)
\end{array}\right.\,, \phantom{AAAAAAA}(y,x)\in\real^n\times\torus^n\,,\nonumber
\\
&&\ts \ham(y,x;\e):=\frac12 |y|^2 +\e f(x)\,,\phantom{AAA}  n\ge 2\,,\ 0<\e<1\,.
\eeqa
As usual, `dot' denotes derivative with respect to `time' $t\in\real$;  $\ham_y$ and  $\ham_x$ denote the gradients with respect to $y$ and  $x$; $|y|^2:=y\cdot y:=\sum_j|y_j|^2$; $\torus^n$ denotes the standard flat torus $\real^n/(2\pi \integer^n)$, and 
the phase space   $\R^n\times \torus^n$ is endowed with the standard symplectic form $dy\wedge dx=\sum_j dy_j\wedge dx_j$.

\nl
In this section, we discuss the  high order normal forms of generic natural systems, especially in neighbourhoods of simple resonances. 

\nl
As standard in perturbation theory, we consider a bounded phase space $\cM\subseteq \R^n\times\T^n$.
By  translating actions and rescaling the parameter $\e$, it is not restrictive to take
\begin{equation}\label{bada}
\cM:={\rm B}\times \T^n\,,\qquad {\rm with}\qquad \DD:=B_1(0):=\{y\in \R^n\st  |y| <1\}\,.
\end{equation} 

\nl
The first step in averaging theory is to construct suitable coverings so as to control resonances where small divisors appear.
Let us recall that {\sl a resonance} ${\cal R}_k$
(with respect to the free Hamiltonian $\frac12 |y|^2$) is the hyperplane
 $\{y\in \real^n: y\cdot k=0\}$, where $k\in\gen$, and its order is given by $\noruno{k}$;
a {\sl double resonance} 
${\cal R}_{k,\ell}$  is the intersection of two   resonances:  ${\cal R}_{k,\ell}={\cal R}_k\cap {\cal R_\ell}$ with 
$k\neq\ell$ in
 $\gen$; the order of ${\cal R}_{k,\ell}$ is given by $\max\{\noruno{k},\noruno{\ell}\}$.

\subsubsection*{Notations}
The real or complex (open) balls of radius $r>0$ and center $y_0\in \R$ or $z_0\in \complex^n$ are denoted by
\beq{palle}
B_r(y_0):= \{y\in\R^n: |y-y_0|<r\}\,,\qquad D_r(z_0):= \{z\in\complex^n: |z-z_0|<r\}\,;
\eeq
if $V\subset \real^n$ and $r>0$, $V_r$  denotes the complex neighborhood of $V$ given by\footnote{$\dst |u|:=\sqrt{u\cdot \bar u}$ denotes the standard Euclidean norm on vectors $u\in\complex^n$ (and its subspaces); `bar', as usual,  denotes complex--conjugated.}
\beq{dico}
V_r := \bigcup_{y\in D}\ D_r(y)\,.
\eeq
We shall also use the notation $\Re(V_r)$ to denote the {\sl real} $r$--neighbourhood of $V\subset \real^n$, namely,
\beq{dire}
\Re(V _r) := V_r\cap \real^n= \bigcup_{y\in V}\ B_r(y)\,.
\eeq
For  a set $V\subseteq \R^n$ and for $r,s>0$, given a function  
$f:(y,x) \in V_r\times \torus^n_s\to f(y,x)$,
we denote
\beq{norme}
{\modulo}f{\modulo}_{V,r,s} = {\modulo}f{\modulo}_{r,s}:=\sup_{V_r\times \T^n_s}|f|\,,
\qquad
\norma f\norma_{V,r,s}=\norma f\norma_{r,s}:=
\sup_{y\in V_r}\sum_{k\in\Z^n} 
|f_k(y)| e^{|k|_{{}_1}s}\,,
\eeq
where $f_k(y)$ denotes the $k$--th Fourier coefficient of   $x\in\torus^n\mapsto f(y,x)$; 
for a function depending only on $y\in V_r$ we set ${\modulo}f{\modulo}_{V,r}={\modulo}f{\modulo}_{r}:=
\sup_{V_r}|f|$.

\subsection{Non--resonant and simply--resonant sets}

Denote  by $\pk$ and $\pko$ the orhogonal projections
\beq{orto}
\pk y:=(y\cdot e_k)\, e_k\,,\qquad \pko y:=y-\pk y\,,\qquad e_k:=k/|k|\,,
\eeq   
and, for any $\K\geq \KO\geq 2$ and  $\a>0$,  define the following sets:
\beqa{neva} 
&&\Rz:=\{y\in \DD:  \min_{k\in \genKO}|y\cdot k|>\ts \frac\a2  \}\,, \ 
\\
\label{sonnosonnoBIS}
&&
\left\{\begin{array}{l}
\Ruk:=
\big\{y\in \DD:
 |y\cdot k|<\ts\a;\,  |\pko y\cdot \ell|> \frac{3 \a \K}{|k|}, \forall
\ell\in {\mathcal G^n_{\K}}\bks\Z k\big\}\,,\quad
(k\in\genKO);
\\
\Ru:=\bigcup_{k\in  \genKO} 	\Ruk\,;
\end{array}\right.
\eeqa
where, as above, $\DD=B_1(0)$.

\nl
Eq.   \equ{neva} implies that $\Rz$ is a $(\a/2)$--non--resonant set up to order $\KO$, i.e.,
\beq{ovvio}
|y\cdot k|>\frac\a2\,,\qquad \forall \ y\in \Rz\,,\ \forall \ 0<|k|\le \KO\,.
\eeq
Indeed,  fix $y\in \Rz$ and $k\in\integer^n$ with $0<|k|\le \KO$. Then,  there exists $\bar k\in \genKO$ and $j\in\integer\bks\{0\}$ such that $k=j\bar k$, so that
$$
|y\cdot k|=|j|\ |\bar k\cdot y|\ge |\bar k\cdot y|>\a/2\,.
$$
From \equ{sonnosonnoBIS} it follows that $\Ruk$
is $(2 \a \K/|k|)$--non resonant modulo $\integer k$ up to order $\K$, namely:
\begin{equation}\label{cipollotto2}
|y\cdot \ell |\ge 2\a\K/|k|\,,
\ \ \ 
\forall k \in \genKO\,,\ 
\forall  y\in \Ruk\,,\ 
\forall \ell\notin \Z k\,, \     |\ell|\leq \K\ .
\end{equation}
Indeed,  fix  $y\in \Ruk$, $k\in\genKO$, $\ell\notin \Z k$ with  $|\ell|\le \K$.
Then,    there exist $j\in\Z\setminus\{ 0\}$
and $\bar\ell\in\genK$ such that $\ell=j\bar\ell$. Hence,   
\beqano
|y\cdot \ell|&=&|j|\, |y\cdot \bar\ell| \ge |y\cdot \bar\ell| =| \pko y \cdot\bar\ell+ \pk y\cdot \bar\ell|
\ge |\pko y\cdot \bar\ell|- \frac{\a \K}{|k|}\\
&>& \frac{3 \a \K}{|k|}  - \frac{\a \K}{|k|} = \frac{2 \a \K}{|k|}\ .
\eeqano
Relations \equ{ovvio} and \equ{cipollotto2}  yield  quantitative control on the small divisors that appear in 
perturbation theory allowing for {\sl high} order averaging theory as we now proceed to show.

\subsubsection*{Averaging}
To perform averaging,  
we need to introduce a few parameters (Fourier cut--offs, a small divisor threshold,  radii of analyticity) and some notation.

\nl
Let 
\beq{dublino}
\left\{\begin{array}{l}
\K\ge  6 \KO\ge 12\,,\quad
 \nu:= \ts\frac92n+2\,, 
\quad
\a:= \sqrt\e \K^\nu\,,
\quad
r_{\rm o}:=\frac{\a}{16 \KO}\,,   
\quad r_{\rm o}':= \frac{r_{\rm o}}2\,,
\\
\ts s_{\rm o}:=s\big(1-\frac1{\KO}\big)\,, 
\  s_{\rm o}':=s_{\rm o}\big(1-\frac1{\KO}\big)\,,
\ 
\ts  s_{\varstar}:=s\big(1-\frac1{\K}\big)\,,
\  s_{\varstar}':=s_{\varstar}\big(1-\frac1{\K}\big)\,,\\
 r_k:={\a}/{ |k|}=\sqrt\e \K^\nu/|k|\,,\quad r_k':=\frac{r_k}2\,,\ s'_k:=|k|_{{}_1}s_{\varstar}'\,,\qquad (\forall\ k\in\genKO)\,.
\end{array}\right.
\eeq
\rem
(i) The action space $\DD$ can be trivially covered by three sets as follows
$$
\DD=\Rz\cup\Ru \cup \Rd\,,\qquad \Rd:=\DD\bks (\Rz\cup\Ru)\,.
$$
As just pointed out, on the set $(\Rz\cup\Ru)\times \T^n$ one can construct detailed, high order normal forms,
while $\Rd$ is a {\sl small set of measure of order $\e^2 \K^\g$} (compare \equ{teheran4} below).

\nl
(ii) It is important to notice that $\Rd$, which is a neighborhood of double resonances of order $\K$,  is a {\sl non perturbative set}, as pointed out in \cite{AKN}. Indeed, consider for simplicity the case $n=2$, 
where the only double resonance is the origin $y=0$. 
Rescaling  variables and time by setting  $y =\sqrt\e {\rm y}$, ${\rm x}=x$,
${\rm t}=\sqrt{\e}t$, the Hamiltonian $t$--flow of $\frac12 |y|^2+\e f(x)$ on $\{y: |y|<\e\}\times \T^2\subseteq \Rd\times \T^2$  is equivalent to the  ${\rm t}$--flow on  $\{|{\rm y}|<1\} \times \T^2$ of the Hamiltonian $\frac12 {\rm y}^2+f({\rm x})$, which {\sl does not depend upon $\e$}. 
\erem

\nl
Next  result is based on 
 `refined Averaging Theory' as presented  in \cite{BCnonlin}. The main technical point in this approach is  the minimal loss of regularity in the angle analyticity domain and the usage of two Fourier cut--offs; for a discussion on these fine points, we refer to the Introduction in \cite{BCnonlin}.

\begin{lemma}[Averaging Lemma] 
\label{averaging} Let $\ham$ be as in \equ{ham} with  $f\in\Bns$ and let  \equ{dublino} hold. There exists a constant
$\bco=\bco(n,s)>1$ such that if $\KO\ge \bco$ the following holds. 

\nl
{\rm (a)}  There exists a real analytic symplectic map
\begin{equation}\label{trota}
\Psi_{\rm o}: \Rz_{r_{\rm o}'}\times \T^n_{s_{\rm o}'} \to 
\Rz_{r_{\rm o}}\times \T^n_{s_{\rm o}} 
\,,
\end{equation}
such that, denoting by $\langle \cdot \rangle$  the 
average over angles $x$,
\beq{prurito}
\hamo(y,x) := \big(\ham\circ\Psi_{\rm o}\big)(y,x)
=\frac{|y|^2}2+\e\big( g^{\rm o}(y) +
 f^{\rm o}(y,x)\big)\,,\quad
\langle f^{\rm o}\rangle=0\,,
\eeq
where  $g^{\rm o}$ and $f^{\rm o}$ are real analytic on $\Rz_{r_{\rm o}'}\times \T^n_{s_{\rm o}'}$ and satisfy 
\beq{552}
| g^{\rm o}|_{r_{\rm o}'}
\leq
\tettao:=
\frac{1}{\K^{6n+1}}\,, \qquad  
\norma f^{\rm o} \norma_{r_{\rm o}',s_{\rm o}'} 
\leq  e^{-\KO s/3}\,.
\eeq
{\rm (b)} For each $k\in \genKO$, there exists a real analytic symplectic map
\begin{equation}\label{canarino}
\Psi_k: 
\Ruk_{r_k'} \times \T^n_{s_{\varstar}'} 
\to 
\Ruk_{r_k} \times \T^n_{s_\varstar} 
\,,
\end{equation}
such that
\beqa{hamk}
\hamk(y,x) &:=& \big(\ham\circ\Psi_k\big)(y,x)\\
&=&\frac{|y|^2}2+\e \big( g^k_{\rm o}(y)+ 
g^k(y,k\cdot x) +
f^k (y,x)\big)\,,\qquad\fproj f^k=0\,,
\nonumber
\eeqa
where
$g^k_{\rm o}$ is real analytic on $\Ruk_{r_k'}$;
 $g^k(y,\cdot)\in\hol_{s'_k}^1$
for every $y\in \Ruk_{r_k'}$ (in particular, $\langle g^k(y,\cdot)\rangle=0$); $f^k $ is real analytic on $\Ruk_{r_k'} \times \T^n_{s_{\varstar}'} $, and:
\begin{equation}\label{cristina}
|g^k_{\rm o}|_{r_k'}
\leq \tettao\,,\qquad
\norma  g^k-\fproj f\norma _{r_k',s'_k} 
\leq  \tettao\,,\qquad 
\norma f^{k} \norma _{r_k',\frac{s_{\varstar}}2} \le
 e^{- \K s/3}\,.
\end{equation}
{\rm (c)} Finally, denoting by $\pi_y$ and $\pi_x$ the projections onto, respectively,  the action variable $y$ and the angle variable $x$, one has
\begin{equation}\label{dunringill}\ts
|\pi_y\Psi_{\rm o}-y|_{r_{\rm o}',s_{\rm o}'}\leq \frac{r_{\rm o}}{2^7 \KO}\,,\quad
|\pi_y\Psi_k-y|_{r_k', s_{\varstar}'}\leq \frac{r_k}{2^7 \K}
\,,
\end{equation}
and, for every fixed $y$, $\pi_x \Psi_{\rm o}(y,\cdot)$,
and $\pi_x \Psi_k(y,\cdot)$ are  diffeomorphisms on $\T^n$.
\end{lemma}

\proof The statements follow from 
Theorem~6.1 in \cite[p. 3553]{BCnonlin} with obvious notational changes, which we proceed to spell out.
The correspondence of symbols between this paper and \cite{BCnonlin} is the following\footnote{In these identities, the first symbol is the one used here, the second one is that used in \cite{BCnonlin}}:

\beqano
&& 
\Rz=\O^0\,;\ \ \Ruk=\O^{1,k}\,,\ \   \ts \frac{|y|^2}2=h(y)\,; \ \ 
\KO={\mathtt K_{{}_1}}\,,\ \ \K={\mathtt K}_{{}_2}\,,
\\
&&
g^{\rm o}= {\rm g}^{\rm o}\,; \ \ 
f^{\rm o}=f^{\rm o}_{\varstar\varstar}\,;\ \ 
g^k_{\rm o}(y)+g^k(y,k\cdot x)={\rm g}^k(y,x)\,;\ \ f^k=f^k_{\varstar\varstar}\,;
\eeqano
the constants $\bar L$ and $L$ in Definition~2.1 in \cite[p. 3532]{BCnonlin} in the present case are $\bar L=L=1$ (since the frequency map here  is the identity map);  
the projection ${\rm p}_{\!{}_{\integer k}}$ introduced in \cite[p. 3529]{BCnonlin} is different from the projection $\fproj$ defined here, the relation between the two being: $\fproj f(k\cdot x)={\rm p}_{\!{}_{\integer k}}f(x)$; 
finally,  the norm $|\cdot|_{D,r,s}$ in \cite[p. 3534]{BCnonlin} corresponds here to the norm $\norma\cdot\norma_{D,r,s}$, hence:
$$|g^k_{\rm o}|_{r_k'}+\norma  g^k-\fproj f\norma _{r_k',s'_k} = \norma {\rm g}^k - {\rm p}_{\!{}_{\integer k}} f\norma_{D^{1,k},r_k/2,s_\star}\,.$$
Now,   Assumption~A in \cite[p. 3533]{BCnonlin} holds. Indeed: 
\begin{itemize}
\item[\footnotesize $\bullet$]
 the action--analyticity radii are the same as in \cite{BCnonlin} (compare \equ{dublino} with Eq.~(140) in \cite{BCnonlin});
\item[\footnotesize $\bullet$]
 the angle--analyticity radii defined here are  the same as in 
Eq.s~(144) and (147) in \cite{BCnonlin} (with different names);
\item[\footnotesize $\bullet$]
In \cite{BCnonlin} it is assumed that $\K\ge 3\KO\ge 6$ (see  Eq. (139) in \cite{BCnonlin}), which in view of \equ{dublino}, is satisfied.
Also $\nu$ in \cite{BCnonlin} is assumed to  satisfy $\nu\ge n+2$, which in \equ{dublino} is defined as  $\nu=\frac92n+2$. 

\item[\footnotesize $\bullet$] By taking $\bco$ big enough  condition (143) is satisfied.

\item[\footnotesize $\bullet$] 
finally, 
to meet the smallness condition~(40) in \cite{BCnonlin}, namely $\e\le r^2/\K^{2\nu}$ (where $r$ is the analyticity radius of the unperturbed Hamiltonian, which here is a free parameter), one can take $r=\K^\nu$  so that  condition~(40)
in \cite{BCnonlin} becomes simply $\e\le 1$.
\end{itemize} 
Thus, Theorem~6.1 of \cite{BCnonlin} can be applied, and 
\equ{prurito} and \equ{hamk} are immediately recognized 
as, 
respectively,  Eq.'s  (145) and (148) in \cite{BCnonlin}. 
Since $\bar\vartheta$ and $\vartheta$ in Eq. 141 of \cite{BCnonlin} are of the form $c(n,s) /\K^{7n+1}$, we see that, by taking $\bco$ big enough, they  can be bounded by $\tettao=1/\K^{6n+1}$ in \equ{552}. Analogously, the exponential estimates on the perturbation functions in (146) and (150) of \cite{BCnonlin} are, respectively, of the form $c(n,s)\, \KO^n e^{-\KO s/2}$ and  
$c(n,s)\, \K^n e^{-\K s/2}$, which, by taking $\bco$ big enough, can be bounded, respectively,  by $e^{-\KO s/3}$ and $e^{-\K s/3}$
as claimed. Thus  (a) and (b) are proven. Finally, from  (71) and (69)  in \cite[p. 3541]{BCnonlin} it follows at once
\equ{dunringill}  and the injectivity of the angle maps.
 \qed

\nl
For high Fourier modes,  a more precise and uniform normal form can be achieved\footnote{This lemma should be  compared with 
Theorem 2.1 in \cite{BCnonlin}.}:

\begin{lemma}[Cosine--like Normal forms] \label{coslike}
Let $\ham$ be as in \equ{ham} with  $f\in\Bns$ satisfying \eqref{P1+} and let  \equ{dublino} hold. 
There exists a constant $\bfco=\bfco(n,s,\d)\ge \max\{\Nf\,,\,\bco\}$ such that 
if $\KO\ge \bfco$ then the following holds.
For any 
  $k\in \genKO$  such that
$\noruno{k}\ge \Nf$, then the Hamiltonian $\hamk$ in \equ{hamk} takes the form:
\begin{equation}\label{hamkc} 
\hamk
=
\frac{|y|^2}2 + \e g^k_{\rm o}(y)+
2|f_k|\e\ 
\big[\cos(k\cdot x +\sa_k)+
F^k_{\! \varstar}(k\cdot x)+
g^k_{\! \varstar}(y,k\cdot x)+
f^k_{\! \varstar} (y,x)
 \big]\,,
\end{equation}
where    $\sa_k$ and $F^k_{\! \varstar}$ are as in Proposition~\ref{pollaio}  and:
\begin{equation}\label{cate}
g^k_{\! \varstar}:=\frac{1}{2|f_k|}\, \big(g^k- \fproj f\big)\,,
\qquad
f^k_{\! \varstar} :=\frac{1}{2|f_k|} f^k\,.
\end{equation}
Furthermore,
 $g^k_{\! \varstar}(y,\cdot )\in\hol_1^1$
(for every $y\in \Ruk_{r_k'}$), $\fproj f^k_{\! \varstar}=0$, and  one has:
\beq{martinaTE}
\norma g^k_{\! \varstar}\norma_{r_k',1}\le 
\tetta:=\frac{1}{\K^{5n}}\,,\qquad\quad
\norma f^k_{\! \varstar} \norma _{r_k',\frac{s_\varstar}2} 
\leq
 e^{-\K s/7}\,.
 \eeq
\end{lemma}
Observe that, under the assumptions of Lemma~\ref{coslike}, by \equ{dublino} and \equ{bollettino1} it is
\beq{bollettino3}
\K\ge 6\KO\ge6\Nf\ge 12\ttcs\ge 12\,. 
\eeq

\proof
First of all observe that   the hypotheses of Lemma~\ref{coslike} imply those of Lemma~\ref{averaging} so that the results of Lemma~\ref{averaging} hold.\\
From \equ{cate} it follows that 
$g^k(y,\sa)=\fproj f(\sa)+ 2 |f_k| g^k_\star(y,\sa)$,
which together with \equ{alfacentauri} and \equ{hamk} of
Lemma~\ref{averaging}, implies immediately the relations \equ{hamkc}.
To prove the first estimate in \equ{martinaTE}, we observe that,
since $\noruno{k}\ge \Nf$,  recalling \equ{dublino} and \equ{bollettino3} one has
\beq{bad}\ts
s'_k = 
|k|_{{}_1}  s\, \big(1-\frac1\K\big)^2 > \Nf s\, \frac45> 1\,.
\eeq
Thus, $g^k_{\! \varstar}(y,\cdot)$
 is bounded on a `large' angle--domain of size larger than 1 and has
zero average (since $g^k_{\! \varstar}(y,\cdot)\in\hol_{|k|_{{}_1}s_{\varstar}'}^1$).
Now, recall the smoothing property \equ{lesso} (with $N=1$),
recall that $\KO\le \K/6$, and take $\bfco$ large enough. Then, 
\begin{align*}
\norma g^k_{\! \varstar}\norma_{r_k',1}&\stackrel{\equ{cate}}{:=}
\frac{1}{2|f_k|}\, \norma  g^k- \fproj f\norma_{r_k',1}
\stackrel{\equ{P1+}}\le
\frac{\noruno{k}^n e^{\noruno{k}s}}{2\d}\, \norma  g^k- \fproj f\norma_{r_k',1}
\\
&\stackrel{(\ref{lesso},\ref{bad})}\le
\frac{\noruno{k}^n e^{\noruno{k}s}}{2\d}\, \norma  g^k- \fproj f\norma_{r_k',s'_k} \cdot  e^{-(s'_k-1)}
\stackrel{\equ{cristina}}\le
\frac{\noruno{k}^n e}{2\d}\, \tettao\ e^{\noruno{k}(s-s_\varstar')}\\
&\stackrel{\equ{dublino}}= 
\frac{\noruno{k}^n e}{2\d}\, \tettao\ e^{\frac{\noruno{k}}{\K} s \big(2-\frac1\K\big)}
\stackrel{\equ{552}}\le 
\frac{\KO^n e}{2\d}\, \frac{1}{\K^{6n+1}}\ e^{2s \frac{\KO}{\K}}\le  \frac{1}{\K^{5n}}\stackrel{\equ{martinaTE}}=\tetta\,.
\end{align*}
Furthermore, possibly increasing  $\bfco$, one also has 
\begin{align*}
\norma f^k_{\! \varstar}\norma_{r_k',\frac{s_\varstar}2}&\stackrel{\equ{cate}}=
\frac{1}{2|f_k|}\, \norma  f^k\norma_{r_k',\frac{s_\varstar}2}
\stackrel{ \equ{P1+}}\le
\frac{\noruno{k}^n e^{\noruno{k}s}}{2\d}\, \norma  f^k\norma_{r_k',\frac{s_\varstar}2}
\stackrel{\equ{cristina}}\le
\frac{\noruno{k}^n e^{\noruno{k}s}}{2\d}\, e^{-\frac{\K s}3}
\\
&\le
\frac{\KO^n}{2\d}\ e^{-\K s\big(\frac13-\frac{\KO}\K\big)}
\le
\frac{\K^n}{2\d\cdot 6^n}\ e^{-\K s/6}
\le e^{-\K s/7}\,. \qquad \qedeq
&
\end{align*}

\subsection{Coverings}

As mentioned in the Introduction, the averaging symplectic maps $\Psi_{\rm o}$ and $\Psi_k$ of Lem\-ma~\ref{averaging} may displace boundaries  by $\sqrt\e\K^\nu$ (compare \equ{dublino} and \equ{dunringill}) so  one cannot use the secular Hamiltonians to describe the dynamics all the way up to the boundary of $\DD\times \T^n$.
Such a problem -- which is essential, for example, in achieving the results described at the end of the Introduction -- may be overcome by introduce {\sl a second covering}, as we proceed now to explain.

\nl
Recall the definitions of $\Rz$ and $\Ruk$ in \equ{neva} and \equ{sonnosonnoBIS}; recall {\equ{dublino}, the notation in \equ{dire} and define 
\beq{rocket}
\Rzt:= \Re (\Rz_{r'_{\rm o}/2})\,,\qquad \Rukt:= \Re(\Ruk_{r'_k/2})\,,\phantom{AAAAa}   (k\in\genKO)\,.
\eeq
Then, the following result holds:
\begin{lemma}\label{sedformosa} {\bf (Covering Lemma)} 
\beqa{surge}
&&\Rz\times \T^n \subseteq \Psi_{\rm o}\big(\Rzt\times \T^n\big)\,,\\
&&\label{surge2}
\Ruk\times\T^n\subseteq  \Psi_k\big(\Rukt\times \T^n\big)\,,\qquad \forall k\in \genKO\,,\\
\label{zucchina}
&& 
\label{sonnosonnosonno}
\Rd:=\DD\bks(\Rz\cup\Ru)\subseteq\bigcup_{k\in \genKO} 
\bigcup_{ \sopra{\ell\in \mathcal G^n_{\K}}{\ell\notin  \Z k}}  \Rd_{k,\ell}\,,
\eeqa
where
\beq{defi}
\Rd_{k\ell}:= \big\{y\in \DD:
 |y\cdot k|<\ts\a;\,  |\pko y\cdot \ell|\le \frac{3 \a \K}{|k|}\big\}\,,\qquad (k\in\genKO\,,\ \ell\in \genK\bks\Z k)\,.
\eeq
\end{lemma}

\rem
(i) 
From the definition of $\Rd$ in \equ{sonnosonnosonno} it follows trivially that $\{{\cal R}^i\}$ is a covering of $\DD$
so that  $\DD= \Rz\cup\Ru\cup \Rd$.

\nl
(ii) Notice that from the definition of   of $\Rukt$ in \equ{rocket}, one has that
\beq{rocket2}
\Rukt_{r'_k/2}\subseteq \Ruk_{r'_k}\,.
\eeq

\nl
(iii) Relations \equ{surge} and \equ{surge2} allow to map back the dynamics of the averaged 
Hamiltonians \equ{prurito} and \equ{hamk}  so as to describe the dynamics also {\sl arbitrarily close to the boundary of the starting phase space}. 
\erem

\nl
For the proof of the Covering Lemma we  shall use  the following immediate consequence of the Contraction Lemma\footnote{Recall the definitions in \equ{palle}; as usual $\overline{A}$ denotes the closure of the set $A$.}:

\begin{lemma}\label{DAD}
Fix $y_0\in\R^n$, $r>0$ and let $\phi:D_{2r}(y_0)\to\C^n$ be a real analytic map  satisfying 
\beq{trasloco}
\sup_{D_{2r}(y_0)}|\phi(y)-y|\le M
\eeq
 for some $0<M<r$.
Then,
$y_0\in \phi(\overline{B_r(y_0)})$. 
\end{lemma}
\proof
Let $V_0:=\overline{B_r(0)}$.
Solving the equation $\phi(y)=y_0$
for some $y\in \overline{B_r(y_0)}$ 
is equivalent to solve the fixed point
equation $w=\psi_0(w):=-{\psi}(y_0+w)$
for $w\in V_0$ having set ${\psi}(y):=\phi(y)-y$. By \equ{trasloco} it follows that $\psi_0: V_0 \to V_0$ and 
by the mean value theorem and Cauchy estimates we get that, for every $w,w'\in V_0$,
$$|\psi_0(w)-\psi_0(w')|=|{\psi}(y_0+w)-{\psi}(y_0+w')|
\leq \frac{M}{r} |w-w'|\,,$$
 showing that  $\psi_0$ is a contraction on $V_0$ (since $M/r<1$)  and the claim follows by the standard Contraction Lemma. \qed

\proof {\bf of \equ{surge}} We start by  proving that
\beq{amicamea}
\forall\ (y_0,x)\in \Rz\times\T^n\,,\  \exists!  \ (y,x_0)\in \Rzt\times\T^n\!:\ \ \Psi_{\rm o}(y,x)=(y_0,x_0)\,.
\eeq
Define 
\beq{zoccolo}
M:=\frac{r_{\rm o}}{2^7 \KO}
 \stackrel{\eqref{dublino}}= \frac{\a}{2^{11} \KO^2}
 <
 \frac{\a}{2^{10}\KO^2}=:r<
  \frac{\a}{2^7 \KO}
   \stackrel{\eqref{dublino}}=
 \frac{r_{\rm o}'}4
 \,.
 \eeq
 Fix    $(y_0,x)\in \Rz\times\T^n$ and
 let $\phi(y):=\pi_y\Psi_{\rm o}(y,x)$. Then,  by \equ{zoccolo}, 
$$\sup_{D_{2r}(y_0)}|\phi(y)-y|
\le \sup_{D_{r'_{\rm o}}(y_0)}|\phi(y)-y|
\le |\pi_y\Psi_{\rm o}-y|_{r_{\rm o}',s_{\rm o}'} 
\stackrel{\equ{dunringill}}\le M\,.$$
Thus, 
by
Lemma~\ref{DAD}, since by \equ{zoccolo} $2r<r'_{{\rm o}}/2$, by definition of $\Rzt$,  we have that 
$$y_0\in \pi_y\Psi_{\rm o}\big(\overline{B_r(y_0)}\times \{x\}\big)
\subseteq 
\pi_y\Psi_{\rm o}\big(\Rzt\times \{x\})\,,
$$ which implies that $\Psi_{\rm o}(y,x)=(y_0,x_0)$ with $x_0\in\T^n$ proving
\eqref{amicamea}. Now, observe that the map $ (y_0,x)\in \Rz\times\T^n \mapsto  (y,x_0)\in\Rzt\times \T^n$ in \equ{amicamea}  is nothing else than the diffeomorphism associated to the near--to--identity generating function $y_0\cdot x+\psi_0(y_0,x)$ of the near--to--identity symplectomorphism $\Psi_{\rm o}$. Thus, for each $y_0\in\Rz$, the map $x \in\T^n\mapsto x_0=x+\partial_{y_0} \psi_0(y_0,x)$
is a diffeomorphism of $\T^n$ with inverse given by $x_0\in\T^n\mapsto x=x_0+\chi(y_0,x_0)$ for a suitable (small) real analytic map $\chi$.  Therefore, given $(y_0,x_0)\in\Rz\times \T^n$,  if we take $x= x_0+\chi(y_0,x_0)$ in \equ{amicamea} we obtain that there exist $(y,x)\in \Rzt\times\T^n$ such that $(y_0,x_0)=\Psi_0(y,x)$, proving \equ{surge}. \qed

\proof {\bf of \equ{surge2}}
The argument is completely analogous: Again, we start by proving  that
\beq{amicomeo}
 \forall\ k\in \genKO\,,\ \forall\ (y_0,x)\in \Ruk\times\T^n\,,\  \exists!  \ (y,x_0)\in \Rukt\times\T^n\!:\ \ \Psi_k(y,x)=(y_0,x_0)\,.
\eeq
Fix $k\in\genKO$ and define 
\beq{zoccoli}
M:=\frac{r_k}{2^7 \KO}
 \stackrel{\eqref{dublino}}= \frac{\a}{2^7|k|\, \K}
 <
 \frac{\a}{2^6 |k| \K}=:r<\frac{r_k'}4 \stackrel{\eqref{dublino}}= \frac{\a}{8|k|}
 \,.
 \eeq
Fix $(y_0,x)\in\Ruk\times \T^n$,  and
 let $\phi(y):=\pi_y\Psi_k(y,x)$.
By \equ{zoccoli},   
$$\sup_{D_{2r}(y_0)}|\phi(y)-y|
\le \sup_{D_{r'_k}(y_0)}|\phi(y)-y|
\le |\pi_y\Psi_k-y|_{r_k',s_\varstar} 
\stackrel{\equ{dunringill}}\le M\,.$$
Thus,  by 
Lemma~\ref{DAD} we have 
$$y_0\in \pi_y\Psi_k\big(\overline{B_r(y_0)}\times \{x\}\big)
\subseteq 
\pi_y\Psi_k\big(\Rukt\times \{x\})\,,
$$ which implies that $\Psi_k(y,x)=(y_0,x_0)$ for some $x_0\in\T^n$ proving
\eqref{amicomeo}.
Now, the argument given in the non--resonant case apply also in this case.
\eproof

\proof {\bf of \equ{sonnosonnosonno}} If $y\in  \Rd$ then, since $y\notin \Rz$,  there exists $k\in \genKO$ such that $|y\cdot k|<\a$, in which case,  since $y\notin \Ru$, there exists $\ell\in \genK\bks\Z k$ 
such that $ |\pko y\cdot \ell|\le \frac{3 \a \K}{|k|}$, hence $y\in  \Rd_{k,\ell}$ for some $k\in \genKO$  and  $\ell\in \genK\bks\Z k$. \qed

\nl
Next, we show that 
 the measure of $\Rd$ is proportional to\footnote{A similar result can be found in \cite[p. 3533]{BCnonlin}.} $\a^2$:
\begin{lemma}\label{coverto}\  
There exists a constant $\cdr=\cdr(n)>1$ such that:
\begin{equation}\label{teheran44} 
\meas (\Rd) \le \cdr\,  \a^2\   \K^{2n} \,.
\end{equation}
\end{lemma}

\proof
Let us estimate the measure of $\Rd_{k,\ell}$ in \equ{defi}.
Denote by  $v\in\R^n$  the projection of $y$ onto  
the plane generated by $k$ and $\ell$
(recall that, by hypothesis, $k$ and $\ell$ are not parallel).
Then,   
\begin{equation}\label{soldatino}
|v\cdot k|=|y\cdot k|<\a\,, \qquad |\proiezione_k^\perp v \cdot \ell|
=|\proiezione_k^\perp y \cdot \ell|
 \le 
3\a\K /|k|\,.
\end{equation}
Set
\beq{bacca}
h:=\pk \ell= \ell -\frac{\ell\cdot k}{|k|^2} k\,.
\eeq
Then, $v$ decomposes in a unique way as
$v=a k+ b h$
for suitable $a,b\in\R$.
By \eqref{soldatino},
\beq{goja}
|a|<\frac{\a}{|k|^2}\,,\qquad
|\pk v\cdot\ell|
=|bh\cdot \ell| \le 3\a\K /|k|\,,
\eeq
and, since $ |\ell|^2 |k|^2-(\ell\cdot k)^2$ is a positive integer (recall, that $k$ and $\ell$ are integer vectors not parallel), 
$$
|h\cdot \ell|
\eqby{bacca}
\frac{ |\ell|^2 |k|^2-(\ell\cdot k)^2  }{|k|^2}
\ge \frac1{|k|^2}\,.
$$
Hence, 
\beq{velazquez}
|b|\le 3 \a \K |k|  \,.
\eeq
Then,  write $y\in \Rd_{k,\ell}$ as $y=v+v^\perp$ with 
$v^\perp$ in the orthogonal
complement of the plane generated by $k$ and $\ell$. Since $|v^\perp |\le |y|< 1$  and $v$ lies in the plane spanned by $k$ and $\ell$ inside a rectangle of sizes of length $2\a/|k|^2$ and $6 \a \K |k|$ (compare \equ{goja} and \equ{velazquez})
we find
\beqno\ts
\meas(\Rd_{k,\ell})\le \frac{2\a}{|k|^2}\, (6 \a \K |k|)\ 2^{n-2}=3\cdot 2^n \, \a^2 \frac{\K}{|k|}\,,\quad  \forall 
\left\{\begin{array}{l}
k\in\genKO\,,\\
 \ell\in  \mathcal G^n_{\K}\bks \Z k\,.
\end{array}\right.
\eeqno
Since $\sum_{k\in\genKO}|k|^{-1}\le c\, \KO^{n-1}$ for a suitable $c=c(n)$,  
and $\KO\le \K/6$, \equ{teheran44}  follows. \qed

\rem
In view of  \equ{teheran44} and  \eqref{dublino}, we have
\begin{equation}\label{teheran4} 
\meas (\Rd) 
\le 
\cdr \, \e\, \K^\g 
\,,\qquad
\gamma:=11n+4
\,.
\end{equation}
Thus, if ${\rm V}_{\! n}=\pi^{\frac{n}2}/\G(1+\frac{n}2)$ denotes the volume of the Euclidean unit ball in $\R^n$, we have that
\beq{piccoletto}
\e<\frac{{\rm V}_{\! n}}{\cdr \K^\gamma}\quad\implies \quad \meas (\Rd) <\meas \DD\,.
\eeq
\erem

\subsection{Normal Form Theorem}
In the normal form around simple resonances the `averaged Hamiltonian'   in \equ{hamk} (i.e., the Hamiltonian obtained disregarding the exponentially small term $f^k$) depends on angles through the linear combination   $k\cdot x$, which,
since $k\in\gen$  defines {\sl a new well--defined angle $\ttx_1\in\torus$}. This fact calls for a linear symplectic change of variables:

\begin{lemma}\label{Fiu} 
Let the hypotheses of Lemma~\ref{coslike} hold.
\\
{\rm (i)} For any  $k\in \genKO$  there exists   a matrix
$\hAA\in\Z^{(n-1)\times n}$
such that\footnote{${\rm SL}(n,\Z)$ denotes the group of $n\times n$ matrices with entries in $\integer$ and determinant 1;
$|M|_{{}_\infty}$, with $M$ matrix (or vector), denotes the maximum norm $\max_{ij}|M_{ij}|$ (or $\max_i |M_i|$).}
\beq{scimmia}
\begin{array}l
\dst \AA:=\binom{k}{\hAA}
=\binom{k_1\cdots k_n}{\hAA}\in\ \ {\rm SL}(n,\Z)\,,\\
|\hAA|_{{}_\infty}\leq |k|_{{}_\infty}\,,\ \ 
|\AA|_{{}_\infty}=|k|_{{}_\infty}\,,\ \ 
|\AA^{-1}|_{{}_\infty}\leq 
(n-1)^{\frac{n-1}2} |k|_{{}_\infty}^{n-1}\,.\phantom{\dst\int}
\end{array}
\eeq
{\rm (ii)} Let $\Fio$ be the  linear, symplectic map on $\R^n\times\torus^n$ onto itself defined by
\begin{equation}\label{talktothewind}
\Fio: (\tty,\ttx) \mapsto (y,x)=
(\AA^T\tty,  \AA^{-1} \ttx) \,.
\end{equation}
Then, 
\beq{azz}
\ttx_1=k\cdot x\,,\qquad\quad  y=\tty_1 k+\hAA^T \hat \tty\,,\phantom{AAAAAAA}  \big[\hat \tty:=(\tty_2,...,\tty_{n})\big]\,.
\eeq
Furthermore, letting\footnote{$\Rukt$ is  defined in \equ{rocket}; recall, also, \equ{dublino}.}
\begin{equation}\label{formentera}\ts
\DDD^k:= \AA^{-T}\Rukt\,,\quad 
\left\{\begin{array}l
\tilde r_k:=\frac{{r_k}}{ \itcu |k|}\\
\tilde s_k:=
\frac{s}{\itcu |k|^{n-1}}
\end{array}\right.
\,,\quad \itcu:=5n(n-1)^{\frac{n-1}2}\,,
\end{equation}
with  $\AA$  as in {\rm (i)}, we find
\begin{equation}\label{fare2}
\Fio: \DDD^k_{\tilde r_k}\times \T^n_{\tilde s_k} 
\to \Rukt_{r_k'/2}\times \T^n_{s_\varstar/2}\,, \qquad
\Fio(\DDD^k\times \torus^n)=\Rukt\times\torus^n\,.
\end{equation}
{\rm (iii)}   $\hamk$ in \equ{hamk}, in the  symplectic variables $(\tty,\ttx)=\big((\tty_1,\hat\tty),\ttx\big)$, takes the form:
\beq{dopomedia}
\cH_k(\tty,\ttx):=\hamk\circ \Fio(\tty,\ttx)=\hamsec_k(\tty,\ttx_1)+ \e
\bar f^k(\tty,\ttx) \,,\quad (\tty,\ttx) \in \DDD^k_{\tilde r_k}\times \T^n_{\tilde s_k}\,,
 \eeq
 where the `secular Hamiltonian' 
 \beq{hamsec}
\hamsec_k(\tty,\ttx_1):= \frac12 |\AA^T\tty|^2+\e g^k_{\rm o}(\AA^T\tty)+
\e g^k(\AA^T\tty,\ttx_1)\,,\quad
\bar f^k(\tty,\ttx):=f^k (\AA^T\tty,\AA^{-1} \ttx)
\eeq
is a real analytic function  for $\tty\in  \DDD^k_{\tilde r_k}$ and\footnote{Recall \equ{dublino}.} $\ttx_1\in \torus_{s'_k}$.
\end{lemma}

\rem\label{ariazz} In the above Lemma~\ref{Fiu} (and also often  in what follows), to simplify symbols, we  may omit the dependence upon $k$ in the notation, but of course $\AA$, $\hAA$ and $\Fio$ {\sl do depend upon the simple resonance label $k\in \genKO$}. 
\erem

\proof {\bf of Lemma \ref{Fiu}}  (i) From B\'ezout's lemma it follows that\footnote{See Appendix~A of \cite[p. 3564]{BCnonlin} for a detailed proof.}:

\nl
{\sl 
given $k\in\integer^n$, $k\neq 0$ there exists a matrix $\AA=(\AA_{ij})_{1\le i,j\le n}$ with integer entries such that $A_{nj}=k_j$ $\forall$ $1\le j\le n$, $\det\AA={\rm gcd}(k_1,...,k_1)$, and 
$|\AA|_{{}_\infty}=|k|_{{}_\infty}$.}

\nl
Hence, 
since $k\in \gen$,   ${\rm gcd}(k_1,...,k_1)=1$,  and \equ{scimmia} follows\footnote{Notice that the bound on $|\AA^{-1}|_{{}_\infty}$ follows from D'Alembert expansion of determinants, observing that
for  any $m\times m$ matrix ${\rm M}$, one   has 
$|\det {\rm M}|\leq m^{m/2} |{\rm M}|_{{}_\infty}^m$}.

\nl
(ii)  $\Fio$ is symplectic since it is generated by  the generating function $\tty\cdot \AA x$. \\
The relations in \equ{azz} follow at once from the  definition of $\Fio$.
\\
Let us prove \equ{fare2}: $\tty\in\DDD^k_{\tilde r_k}$ if and only if $\tty=\tty_0+z$ with $\tty_0\in\DDD^k$ and $|z|<\tilde r_k$. Thus, 
 $$|\AA^Tz|\stackrel{\equ{scimmia}}\le n |k|  |z|<n  |k| \tilde r_k\stackrel{ \equ{formentera}}< \frac{r_k}4
 \stackrel{ \equ{dublino}}=\frac{r'_k}2\,.
 $$ 
Since, by definition of $\DDD^k$, $\AA^T \tty_0\in \Rukt$, we have that 
 $\AA^T\tty\in \Rukt_{r_k'/2}$.
 \\
Let, now, $\ttx$ belong to $\T^n_{\tilde s_k}$.Then, for any $1\le j\le n$, recalling the definitions of $s_\varstar$ and $s_\varstar'$ in \equ{dublino}, we find 
$$
\Big|\Im (\AA^{-1} \ttx)_j\Big|= \Big| \sum_{i=1}^n (\AA^{-1})_{ij} \Im \ttx_j\Big|\stackrel{\equ{scimmia}}{<} 
n(n-1)^{\frac{n-1}2} |k|^{n-1} \tilde s_k \stackrel{\equ{formentera}}{\le }\frac{s_\varstar}{2}< s_\varstar'\,.
$$
Thus, $\AA^{-1} \ttx$ belong to $\T^n_{s_\varstar'}$, and \equ{fare2} follows. 

\nl
(iii) Eq.'s \equ{dopomedia}--\equ{hamsec} follow immediately from the definition of the symplectic map $\Fio$ in  \equ{talktothewind} and \equ{azz}. The statement on the angle--analyticity domain of $\hamsec_k$ follows from part (b) of Lemma~\ref{averaging}.
\qed 

\nl
We summarize the above lemmata in the  following

\begin{theorem}[Normal Form Theorem]
\label{normalform}
Let $\ham$ be as in \equ{ham} with  $f\in\Bns$ satisfying \eqref{P1+} with $\Nf$ as in \equ{enne}, and let  \equ{dublino} hold. 
There exists a constant\footnote{$\bco$ is defined in Lemma~\ref{averaging}.} $\bfco=\bfco(n,s,\d)\ge \max\{\Nf\,,\,\bco\}$ such that, 
if $\KO\ge \bfco$, $k\in\genKO$,  and  $\DDD^k$, $\tilde r_k$, $\tilde s_k$ are as in \equ{formentera}, then
there exist
 real analytic symplectic maps
\beq{trota2}
\Psi_{\rm o}: \Rz_{r_{\rm o}'}\times \T^n_{s_{\rm o}'} \to 
\Rz_{r_{\rm o}}\times \T^n_{s_{\rm o}} 
\,,
\qquad
\Psi^k:
\DDD^k_{\tilde r_k}\times \T^n_{\tilde s_k} 
\to 
\Ruk_{r_k} \times \T^n_{s_\varstar}
\eeq
having the following properties.

\nl
{\rm (i)} 
$
\hamo(y,x) := \big(\ham\circ\Psi_{\rm o}\big)(y,x)
=\frac{|y|^2}2+\e\big( g^{\rm o}(y) +
 f^{\rm o}(y,x)\big)$,
with  $g^{\rm o}$ and $f^{\rm o}$ satisfying
\equ{552} and  $\langle f^{\rm o}\rangle=0
$.

\nl
{\rm (ii)} 
\beq{colosseum3cippa}
\cH_k(\tty,\ttx):=\ham\circ \Psi^k(\tty,\ttx)=\hamsec_k(\tty,\ttx_1)+ \e
\bar f^k(\tty,\ttx) \,,\quad (\tty,\ttx) \in 
\DDD^k_{\tilde r_k}\times \T^n_{\tilde s_k}\,,
 \eeq
 where
 \beq{hamseccippa}
\hamsec_k(\tty,\ttx_1):= \frac12 |\AA^T\tty|^2+\e \mathtt g^k_{\rm o}(\tty)+
\e \mathtt g^k(\tty,\ttx_1)
\eeq
is a real analytic function  for $\tty\in  \DDD^k_{\tilde r_k}$ and $\ttx_1\in \torus_{s'_k}$.
In particular
 $\mathtt g^k(y,\cdot)\in\hol_{s'_k}^1$
for every $y\in \DDD^k_{\tilde r_k}$. Furthermore, the following estimates hold:
\begin{equation}\label{cristinacippa}
|\mathtt g^k_{\rm o}|_{\tilde r_k}
\leq \tettao=\frac{1}{\K^{6n+1}}\,,\qquad
\norma  \mathtt g^k-\fproj f\norma _{{\tilde r_k},s'_k} 
\leq  \tettao\,,\qquad 
\norma \bar f^{k} \norma _{{\tilde r_k},{\tilde s_k}} \le
 e^{- \K s/3}\,.
\end{equation}
{\rm (iii)} If  $   \noruno{k}\geq \Nf$, there exists  $\sa_k\in[0,2\pi)$ such that 
\begin{equation}\label{hamkccippa} 
\cH_k
=
 \frac12 |\AA^T\tty|^2+\e \mathtt g^k_{\rm o}(\tty)+
2|f_k|\e\ 
\big[\cos(\ttx_1 +\sa_k)+
F^k_{\! \varstar}(\ttx_1)+
\mathtt g^k_{\! \varstar}(\tty,\ttx_1)+
\mathtt f^k_{\! \varstar} (\tty,\ttx)
 \big]\,,
\end{equation}
where   $F^k_{\! \varstar}$ is as in Proposition~\ref{pollaio} and satisfies
$F^k_{\! \varstar}\in\hol_1^1$ and
$
\modulo F^k_{\! \varstar} \modulo_1\leq 2^{-40}$.\\
Moreover,
 $\mathtt g^k_{\! \varstar}(y,\cdot )\in\hol_1^1$
(for every $y\in \DDD^k_{\tilde r_k}$), $\fproj\mathtt f^k_{\! \varstar}=0$, and  one has 
\beq{martinaTEcippa}\ts
\norma \mathtt g^k_{\! \varstar}\norma_{\tilde r_k,1}\le 
\tetta=\frac{1}{\K^{5n}}
\,,\quad\qquad
\norma \mathtt f^k_{\! \varstar} \norma _{\tilde r_k,\tilde s_k} 
\leq
e^{-\K s/7}\,.
\eeq
\end{theorem}
\proof 
The first relation in \equ{trota2} is \equ{trota}. Define
\beq{pippala}
\Psi^k:= \Psi_k\circ \Fio\,.
\eeq
Then,  since $s_\varstar/2<s_\varstar'$ (compare \equ{dublino}), by
 \equ{fare2}, \equ{rocket2} we get the second relation in \equ{trota2}.

\nl
{\rm (i)}
follows from point {\rm (a)}  of Lemma~\ref{averaging}.

\nl
{\rm (ii)}
\equ{colosseum3cippa},
\equ{hamseccippa} and \equ{cristinacippa} follow from, respectively,
\equ{dopomedia},
\equ{hamsec}, \equ{cristina}
and point (ii) of Lemma~\ref{Fiu} setting
\beq{ggg}
\mathtt g^k_{\rm o}(\tty):=
g^k_{\rm o}(\AA^T\tty)\,,\qquad
\mathtt g^k(\tty,\ttx_1)
:= g^k(\AA^T\tty,\ttx_1)\,.
\eeq
{\rm (iii)}
follows by Proposition~\ref{pollaio} and Lemma~\ref{coslike}.
In particular 
\equ{hamkccippa} follows from \equ{hamkc}.
Furthermore,
\begin{equation}\label{catecippa}
\mathtt g^k_{\! \varstar}:=\frac{1}{2|f_k|}\, \big(\mathtt g^k- \fproj f\big)\,,
\qquad
\mathtt f^k_{\! \varstar} :=\frac{1}{2|f_k|} \bar f^k
\end{equation}
and
noting that 
$\mathtt g^k_{\! \varstar}(\tty,\ttx_1)
= g^k_{\! \varstar}(\AA^T\tty,\ttx_1)$ and that, by \equ{hamsec}, 
$\mathtt f^k_{\! \varstar}(\tty,\ttx)=f^k_{\! \varstar}(\AA^T\tty,\AA^{-1} \ttx)
$, we see that
\equ{martinaTEcippa} follows from \equ{martinaTE} and \equ{fare2}.
\qed

\section{Generic Standard Form at simple resonances}

In this final section we show that {\sl  the secular Hamiltonians $\hamsec_k$  \equ{hamsec} in Theorem~\ref{normalform}
can be symplectically put into a suitable standard form, uniformly in  $k\in\genKO$}

\nl
The precise definition of `standard form' is taken from \cite{BCaa23}, where the analytic properties 
of   action--angle variables of such  Hamiltonian systems are discussed.

\begin{definition}\label{morso}
Let $\hat D \subseteq \R^{n-1}$ be a bounded domain,  $\Ro>0$ and $D:=  (-\Ro ,\Ro ) \times\hat  D
$. We say that the real analytic Hamiltonian $\Hpend$ is in Generic Standard Form with  respect to the symplectic 
variables $(p_1,q_1)\in (-\Ro ,\Ro )\times\torus$ and  `external actions' 
$$\hat p=(p_2,...,p_n)\in \hat D$$ if $\Hpend$ has the form
 \beq{pasqua}
  \Hpend(p,q_1)=
\big(1+ \cin(p,q_1)\big) p_1^2  
  +\Gm(\hat p, q_1)
  \,,
\eeq
where:

\begin{itemize}

\item[\bolla]  $\cin$ and $ \Gm$ are real analytic functions defined on, respectively, $D_\ro\times\T_\so$ and $\hat D_\ro\times \T_\so$ for some $0<\ro\leq\Ro$ and $\so>0$;

\item[\bolla]  
$\Gm$ has zero average and there exists a 
 function $\GO$ (the `reference potential')  depending only on $q_1$ such that , for some\footnote{Recall Definition~\ref{buda}.} $\morse>0$, 
 \begin{equation}\label{A2bis}
\GO\  \  \mbox{is} \ \
\morse {\rm \text{--}Morse}\,,\qquad \langle \GO\rangle=0\,;
\end{equation}
\item[\bolla] 
 the following estimates hold:
 \beq{cimabue}
 \left\{\begin{array}{l} \dst \sup_{\torus^1_\so}|\GO|\le \suca\,,\\
 \dst \sup_{\hat D_\ro\times \torus^1_\so}|\Gm-\GO| \leq
\suca
 \lalla
\,,\quad{\rm for\ some}\quad 0<\suca\le \ro^2/2^{16} 
\,,\ \ 0\le \lalla<1\,,
\\
\dst \sup_{D_\ro\times \torus^1_\so}|\cin| \leq
\lalla\,.
\end{array}\right.
\eeq
\end{itemize}
\end{definition}
We shall call $(\hat  D,\Ro,\ro,\so,\morse,\suca,\lalla)$ {\sl the analyticity characteristics of $\Hpend$
with respect to the unperturbed potential $\GO$}.

\rem\label{trivia}  
 If $\Hpend$ is in Generic Standard Form, then the parameters $\morse$ and $\suca$ satisfy the relation\footnote{By \equ{cimabue}, $\morse \le |\bar \Gm(\sa_i)- \bar \Gm(\sa_i)|\le 2 \max_\torus|\bar \Gm|\le 2\suca$.}
\beq{sucamorse}\ts
\frac\suca\morse\ge \frac12\,.
\eeq 
Furthermore,  one can always
fix  $\varpi\geq 4$  such that:
  \begin{equation}\label{alce}\ts
\frac{1}{\varpi}\leq \so\leq 1\,,\qquad
1\leq
\frac{\Ro}{\ro}\leq \varpi\,,\qquad
\frac{1}{2}\leq
 \frac{\suca}{\morse }
\leq \varpi \,.
\end{equation}
Such a  parameter $\varpi$ rules the main scaling properties of these Hamiltonians. 
\erem

\subsection{Main theorem}

\nl
In the following we shall often use the  following notation: 
If $w$ is a vector with $n$ or $2n$ components, $\hat w=(w)^{\widehat{}}$ denotes the last $(n-1)$ components; 
if $w$ is vector with $2n$ components, $\check w=(w)^{{\!\!\widecheck{\phantom{a}}}}$ denotes the first $n+1$ components. 
Explicitly:
\beq{checheche}
w=(y,x)=\big((y_1,...,y_n),(x_1,...x_n)\big)\quad \Longrightarrow\quad
\left\{
\begin{array}{l}
\hat w=(w)^{\widehat{}}=(x_2,...,x_n)=\hat x\,,\\
\hat y=\,(y)^{\widehat{}}=(y_2,...,y_n)\,,\\
\check w=(w)^{{\!\!\widecheck{\phantom{a}}}}=(y,x_1)\,,\\
w=(\check w,\hat w)\,.
\end{array}
\right.\eeq

\dfn{dadaumpa}
Given a domain $\hat {\rm D}\subseteq \R^{n-1}$,
we denote by
 $\Gdag$   the 
abelian group of 
symplectic diffeomorphisms $\Psi_{\! \ta}$ 
of $(\R\times\hat {\rm D})\times \R^n$
given by
\beqno
(p,q)\in(\R\times\hat {\rm D})\times \R^n\stackrel{\Psi_{\! \ta}}\mapsto (P,Q)=
(p_1+\ta(\hat p),\hat q,q_1,\hat q-q_1\partial_{\hat p}
\ta(\hat p))\in\real^{2n}\,,
\eeqno
with $\ta:\hat {\rm D}\to \R$ smooth.
 \edfn

\rem
\label{alice} 
The group properties of $\Gdag$ are trivial: 
$$
{\rm id}_{\Gdag}=\Psi_{\! 0}\,,\qquad   
\Psi_{\! \ta}^{-1}=\Psi_{\! -\ta }\,,\qquad
\Psi_{\! \ta}\circ\Psi_{\! \ta'}=\Psi_{\! \ta+\ta'}\,.
$$
Notice, however,  that, unless $\partial_{\hat p} \ta\in\integer^{n-1}$,
maps in $\Psi_{\! \ta}\in \Gdag$  {\sl do not induce  well defined angle maps}  
 $q\in \torus^n\mapsto (q_1,\hat q-q_1\partial_{\hat p}
\ta(\hat p))\in\torus^n$.
\erem
Now, let $f\in\Gns$ satisfy\footnote{Recall that by Lemma~\ref{telaviv} such $\d$ and $\b$ always exist.}   \eqref{P1+} and \equ{P2+} for some $0<\d\le 1$ and $\b>0$ with $\Nf$ defined in \equ{enne}, let $k\in\genKO$, recall \equ{dublino}  and define the following parameters\footnote{Here and in what follows we shall  not always indicate explicitly the dependence upon $k$.
Recall the definitions of $\itcu$, $\hAA$ and $\ttcs$ in, respectively, \equ{formentera}
 Lemma~\ref{Fiu} and \equ{bollettino1}.}
\beq{cerbiatta}
\begin{array}l
\Ro={\a}/{|k|^2}={\sqrt\e \K^\nu}/{|k|^2}\,,\quad \itcd=4\,  n^{\frac32} \itcu\,, 
\quad \ro= {\Ro}/{\itcd}\,,
\quad\e_k=\frac{2\e}{|k|^2}\,,
\phantom{\dst\int}
\\
\hat D =
\big\{ \hat\act\in\R^{n-1}: \ 
|\proiezione_k^\perp \hAA^T \hat\act|<1\,,
\ \dst
\min_{\sopra{\ell\in \genK}{\ell \notin \Z k}}
\big| \big(\proiezione_k^\perp \hAA^T \hat\act\big)\cdot \ell\big|
\geq  {\textstyle \frac{3\a\K}{|k|}}   \big\}\,,\ D=(-\Ro,\Ro)\times \hat D\,,
\\
 \morse=\casitwo{\e_k \b,}{\noruno{k}<\Nf}{\e_k |f_k|,}{\noruno{k}\ge \Nf}\,,\,\qquad
 \chk=\casitwo
{1 \,,}{\noruno{k}<\Nf}
{ |f_k|\,,}{\noruno{k}\ge \Nf}\,,
\quad \suca=\ts \ttcs \e_k\,  \chk\,,    
\\
\so=\casitwo{\min\{\frac{s}2,1\}\,,}{\noruno{k}<\Nf}{1\,,}{\noruno{k}\ge \Nf}\,,\quad
\ts
\chs:=\casitwo{s'_k\,,}{\noruno{k}<\Nf\,,}{1\,,}{\noruno{k}\ge \Nf}\,, 
\quad \dst
\lalla=\frac{1}{\K^{5n}}\,. 
\end{array}
\eeq

\begin{theorem}[Generic Standard Form at simple resonances]\label{sivori}\  \\
Let $\ham$ be as in \equ{ham} with $f\in\Gns$ satisfying   \eqref{P1+} and \equ{P2+} for some $0<\d\le 1$ and $\b>0$ with $\Nf$ defined in \equ{enne}.
Assume that\footnote{$\bfco$ is defined in Theorem~\ref{normalform}.}  $\KO\ge \max\{\itcd,\bfco\}$. Then,  with the definitions given in \equ{cerbiatta}, the following holds for all $k\in\genKO$.

\nl
{\rm (i)}
There exists a real analytic
 symplectic transformation 
\beq{diamond}
\Phi_\diamond:(\ttp,\ttq)\ \in D\times \R^n \to
(\tty,\ttx)=\Phi_\diamond(\ttp,\ttq)\in \R^{2n}\,,
\eeq
such that: $\Phi_\diamond$ fixes $\hat\ttp$ and\footnote{I.e., in \equ{diamond} it is $\tty=\hat\ttp,\ttx_1=\ttq_1$.} $\ttq_1$; 
for every $\hat\ttp\in\hat D$ the map  $(\ttp_1,\ttq_1)\mapsto (\tty_1,\ttx_1)$ is symplectic; 
the $(n+1)$--dimensional map\footnote{Recall the notation in \equ{checheche}.}  $\check\Phi_\diamond$ depends only on the first $n+1$ coordinates $(\ttp,\ttq_1)$, is $2\pi$--periodic in $\ttq_1$ 
and, if $\DDD^k= \AA^{-T}\Ruk$ and $\hamsec_k$ are as in Theorem~\ref{normalform},  one has\footnote{
$r_k$ and $s'_k$ are defined in \equ{dublino}, $\tilde r_k$ in \equ{formentera}.} 
\beq{tikitaka}
\begin{array}l
 \check\Phi_\diamond: 
 D_{\chr}\times \torus_\chs\ \, \to
\DDD^k_{\tilde r_k}\times \T_\chs\,, {\phantom{\dst \int}}
\\
\hamsec_k\circ \check\Phi_\diamond(\ttp,\ttq)=:
\textstyle{\frac{|k|^2}{2}}( \Hsharp(\ttp,\ttq_1)+
\hzk(\hat\ttp))\,,  {\phantom{\dst \int}}
\\ 
\sup_{\hat\ttp\in \hat D_{2\ro}}
\big|\ts\hzk(\hat\ttp)-
\htk(\hat\ttp)
\big|
\leq
\ts 6\, \e_k\lalla\,,\qquad\  \htk(\hat\ttp):={\ts \frac{1}{|k|^2}} 
\normadue \proiezione^\perp_k \hAA^T \hat\ttp\normadue^2\,.
\end{array}
\eeq
{\rm (ii)} 
$\Hsharp$ in \equ{tikitaka} is  in Generic Universal Form according to Definition~\ref{morso}:
\beqno
\Hsharp(\ttp,\ttq_1)=\big(1+\cins(\ttp,\ttq_1)\big)\,  \ttp_1^2 +  \Gf(\hat\ttp,\ttq_1)\,,
\eeqno
having reference potential 
 \beq{paranoia}
\GO=\bGf:= \e_k\, \fproj f\,,
\eeq  
analyticity characteristics  given in \equ{cerbiatta}, and 
$\upkappa$
verifying \equ{alce}  with
\beq{kappa}\ts
\upkappa=\upkappa(n,s,\b):=\max\big\{\itcd\,, 4\ttcs\,, \ttcs/\b
\big\}\,.
\eeq
{\rm (iii)} Finally, $\Phi_\diamond=\Fiuno\circ\Fidue\circ\Fitre$, where\footnote{Recall Definition~\ref{dadaumpa}.}: $\Fiuno:=\Psi_{\!\giuno}\in\Gdag$ with 
$\giuno(\hat \ttp):=-\ts \frac1{|k|^2}{(\hAA k)\cdot \hat \ttp}$; $\Fitre:=\Psi_{\!\gitre}\in\Gdag$  for a suitable real analytic function $\gitre(\hat\ttp)$ satisfying
\beqno\ts
|\gitre|_{4 \chr}< \frac{\e_k\chk}{\chr}\lalla\,,
\eeqno 
and $\Fidue(\ttp,\ttq)=(\ttp_1+\upeta_{{}_2},\hat \ttp, \ttq_1,\hat \ttq+\upchi_{{}_2})$ for suitable real analytic functions $\upeta_{{}_2}=\upeta_{{}_2}(\hat \ttp,\ttq_1)$ and 
$\upchi_{{}_2}=\upchi_{{}_2}(\hat \ttp,\ttq_1)$ 
 satisfying
\beqno\ts
|\upeta_{{}_2}|_{4 \chr,\chs}< \frac{\e_k\chk}{\chr}\lalla\,,\qquad |\upchi_{{}_2}|_{2 \chr,\chs}< \frac{4\e_k\chk}{\chr^2}\,\lalla
\,.
\eeqno
\end{theorem}

\rem\label{rampulla} (i)  One of the main point of the above theorem is that the parameter $\kappa$ in \equ{kappa} 
{\sl does not depend on $k$}. Incidentally, we point out that $\kappa$  depends (indirectly) also on $\d$, since $\d$ appears  in the definition of $\Nf$ and $\b$ is the uniform Morse constant of the first $\Nf$ reference potentials.

\nl
(ii)
Note that by \equ{cerbiatta},
\equ{tikitaka}, \equ{dublino} 
and \equ{bollettino1}
\begin{equation}\label{fangorn}\ts
\min\big\{\frac{s}2,1\big\}\le \so\leq \chs\le {s'_k}\,.
\end{equation}
In particular, the composition $\hamsec_k\circ \check\Phi_\diamond$ is well defined; compare Theorem~\ref{normalform}--(ii).\\
As for the action analyticity radii, notice that, by the definitions in \equ{dublino}, \equ{formentera} and \equ{cerbiatta}, one has  
\beq{bollettino2}
r_k=\Ro\, |k|\,,\qquad
\tilde r_k= \frac{\Ro}\itcu\,.
\eeq
(iii) The three maps which define  $\Phi_\diamond$ have the following purposes:
The first one is needed to decouple the `kinetic  energy' of the 1--d.o.f.  secular system; the second one is introduced so as to get a purely positional 1--dimensional potential; finally,  the third one  puts the momentum coordinate of the equilibria in 0.

\nl
(iv) The proof is fully constructive and the explicit definition of $\Hsharp$ is given in \equ{cins}, 
\equ{fpe}, \equ{guaito},
\equ{limone}, \equ{pontediferro} and  \equ{hamsecu} 
below. 
\erem

\subsection{Proof of the main theorem}
The proof is articulated  in three lemmata. \\
The first lemma shows how to `block--diagonalize' the kinetic energy. For $k\in \genKO$, recall the definition of the matrices 
 $\AA$ and $\hAA$ in \eqref{scimmia}, and define\footnote{ ${\rm I}_{m}$ denotes the $(m\times m)$--identity matrix and recall the notation in \equ{checheche}.}
\beqa{centocelle}
&&
\tty= {\rm U}\ttY:= 
\left(\begin{matrix}
1  &  - \frac1{|k|^2}(\hAA k)^T \cr 0 & \quad{\rm I}_{{}_{n-1}} \cr
\end{matrix}\right)\ttY
\,,\qquad
{\rm i.e.}\qquad 
\casi{\tty_1=\ttY_1 - \frac1{|k|^2}{\hAA k\cdot \hat \ttY}\,,}
{\hat \tty=\hat \ttY\,.}
\eeqa
Then, one has

\begin{lemma}\label{phi1} {\rm (i)}
 Let $\Fiuno$ be the map $\Fiuno(\ttY,\ttX)=({\rm U}\ttY,{\rm U}^{-T}\ttX)$.
Then, $\Fiuno$ is symplectic and 
\beq{finocchio} \DDD^k = {\rm U} \ZZ \,,\qquad 
\Fiunoc: \ZZ_{4 \chr}\times\T_\chs\to \DDD^k_{\tilde r_k}\times \T_\chs\,.
\eeq
{\rm (ii)}  Let 
\beq{hamsecu}
\left\{
\begin{array}{l}
\Guo:=  {\textstyle \e_k} g^k_{\rm o}(\AA^T {\rm U} \ttY)\,,\quad  
\Gu(\ttY,\ttX_1):= {\textstyle \e_k} 
g^k(\AA^T {\rm U} \ttY,\ttX_1)\,,\\ \ \\
 \hamsecu(\ttY,\ttX_1):=
 \ttY_1^2+ \Guo(\ttY)+\Gu(\ttY,\ttX_1)\,,\qquad \langle \Gu(\ttY,\cdot)\rangle=0\,.
 \end{array}\right.
\eeq
Then, if $\hamsec_k$ is as in \equ{hamsec}, one has   
\beq{spigolak}
\hamsec_k\circ \Fiunoc(\ttY,\ttX_1)= {\textstyle\frac{|k|^2}{2}} \, \hamsecu(\ttY,\ttX_1)+ {\textstyle \frac12}
\normadue \proiezione^\perp_k \hAA^T \hat \ttY\normadue ^2\,,
\eeq 
with $\hamsecu$ real analytic on $\ZZ_{4 \chr}\times\T_\chs$ and $\langle \Gu(\ttY,\cdot)\rangle=0$, and 
the following estimates~hold\footnote{$\tetta$ is defined in \equ{martinaTE}. Notice that, by \equ{cerbiatta},  $\chi_{{}_k}\le 1$ for all $k$.}:
\beq{betta}
|\Guo|_{4 \chr}\le \bettao:=2 \e_k \tetta= \frac{2\e_k}{\K^{5n}}\,,\qquad
| \Gu-\bGf|_{4 \chr,\chs}\le \betta:=\chi_{{}_k}\bettao\le\bettao
\,.
\eeq
\end{lemma}

\proof
(i) $\Fiuno$ is symplectic since it is generated by the generating function ${\rm U}\ttY\cdot \ttx$.\\
From the definitions of $\AA$ and ${\rm U}$ in, respectively, \equ{scimmia} and \equ{centocelle}, it follows  
\beq{perpieta}
(\AA^T{\rm U}) \ttY=\ttY_1 k + \hAA^T\hat \ttY-  {\textstyle \frac{(\hAA k)\cdot \hat \ttY}{|k|^2}k}=
\ttY_1 k + \hAA^T\hat \ttY- {\textstyle \frac{  \hAA^T \hat \ttY \cdot k}{|k|^2}k}
=\ttY_1 k + \proiezione_k^\perp \hAA^T \hat \ttY\,.
\eeq	
Thus, $\tty= (\AA^T{\rm U}) \ttY$ if and only if $\tty\cdot k=\ttY_1 |k|^2$ and $\proiezione_k^\perp \tty=\proiezione_k^\perp
\hAA^T \hat \ttY$, which is equivalent to say  $(\AA^T{\rm U}) \ZZ =\Ruk$, which in view of   \equ{formentera}, is equivalent to
$\DDD^k = {\rm U} \ZZ$. Now,  by~\equ{scimmia},  
\beq{UU}
|{\rm U}|, |{\rm U}^{-1}|\le n^{\frac32}\,,
\eeq 
where, as usual, for a matrix $M$ we denote by
$\dst |M|=\sup_{u\neq 0} |Mu|/|u|$  the standard  operator norm. 
Thus, by   \equ{cerbiatta} and \equ{bollettino2} we have (for complex $z$)
\beq{chitikaka}
|z|<4\ro \ \ \implies
\ \ 
|{\rm U} z|< n^{\frac32} 4 \ro = 4 n^{\frac32}\, \frac{\Ro}{\itcd}= \frac{\Ro}{\itcu}= \tilde r_k\,,
\eeq
which, since $\ttX_1=\ttx_1$,  implies that 
$\Fiunoc: \ZZ_{4 \chr}\times\T_\chs\to \DDD^k_{\tilde r_k}\times \T_\chs$, proving  \equ{finocchio}.

\nl
(ii)  By the previous item, the composition 
 $\hamsec_k\circ \Fiunoc$ is well defined and analytic on $\ZZ_{4 \chr}\times\T_\chs$.
From \equ{perpieta} it follows that
$|\AA^T{\rm U} \ttY|^2=|k|^2  \ttY_1^2  + \normadue \proiezione^\perp_k \hAA^T \hat \ttY\normadue ^2$,
and \equ{spigolak}  follows. Notice that since   $g^k(y,\cdot)\in\hol_{s'_k}^1$ (compare Lemma~\ref{averaging}), $\Gu$ has zero average over $\T$.
\\
By the definition of  $\Guo$ and $\Gu$  in \equ{hamsecu}, by \equ{chitikaka},  \equ{cristina} in\footnote{Recall that, by \equ{dublino}, \equ{formentera},  $\tilde r_k<r'_k=r_k/2$.
Recall also the definitions of  
$\tettao$ and $\tetta$ in \equ{552} and 
\equ{martinaTE}.} Lemma~\ref{averaging}, the estimates 
on $|\Guo|_{4 \chr}$ and on $| \Gu-\bGf|_{4 \chr,\chs}$ for $\noruno{k}<\Nf$  in  \equ{betta} follow.
The estimate for $\noruno{k}\ge\Nf$ in \equ{betta} follows from Lemma~\ref{coslike}: see in particular \equ{cate}, \equ{martinaTE}
and \equ{alfacentauri}. \qed

\nl
Next lemma shows how one can  remove the dependence on $\ttY_1$ in  the potential
$\Gu$.

\begin{lemma}\label{avogado} If let $\K\ge \itcd$ 
then,
\beq{tettapic}\ts
 \frac{\bettao}{\chr^2}
 <\frac1{2^{10}}\  \frac\chs{\pi+\chs}<1
 \,,
\eeq
and  the following statements hold.
\\
{\rm (i)} The fixed point  equation
\beq{fpe}
\pp = -{\textstyle \frac12} \, \partial_{\ttY_1} \Guo(\pp,\hat\ttP) -{\textstyle \frac12} \, \partial_{\ttY_1}\Gu(\pp,\hat\ttP,\ttQ_1)
\eeq
has a unique solution $\pp:(\hat \ttP,\ttQ_1)\in \hat\ZZ\times \torus\mapsto \pp(\hat\ttP,\ttQ_1)\in\real$ real analytic on $\hat\ZZ_{4 \chr}\times \torus_\chs$, satisfying
\beq{pitale}\ts
|\pp|_{4 \chr,\chs}<\frac{\bettao}{3 \chr}\,.
\eeq
Furthermore, if we define 
\beq{guaito}
\left\{
\begin{array}{l}
\pp_{\rm o}(\hat \ttP):=\langle \pp(\hat\ttP,\cdot)\rangle
\\
\tilde\pp:=\pp-\pp_{\rm o}
\end{array}\right.\,,\qquad
\left\{
\begin{array}{l}
\dst \phi(\hat \ttP,\ttX_1):= \int_0^{\ttX_1} \tilde\pp(\hat \ttP,\sa)d\sa\\
\hat\qq(\hat\ttP,\ttQ_1):=-\partial_{\hat\ttP} \, \phi(\hat\ttP,\ttQ_1)
\end{array}\right.
\eeq
then, $\ttQ_1\to\hat\qq(\hat\ttP,\ttQ_1)$ is a real analytic periodic function, and
one has
\beq{tess}\ts
|\pp_{\rm o}|_{4 \chr}< \frac13\, \frac{\bettao}{\chr}\,,\qquad\quad
|\tilde\pp|_{4 \chr,\chs}< \frac13\, \frac{\betta}{\chr}\,,\qquad |\hat \qq|_{2 \chr,\chs}< \frac{\betta}{6 \chr^2}\,(\pi+\chs)
\,.
\eeq
{\rm (ii)} 
The real analytic symplectic map $\Fidue$ generated by 
$ \ttP\cdot\ttX+  \phi(\hat\ttP,\ttX_1)$, namely, 
\beq{Fidue}
\Fidue:(\ttP,\ttQ)\mapsto (\ttY,\ttX) \quad{\rm with}\quad
\casi{\ttY_1=\ttP_1+ \tilde \pp(\hat \ttP,\ttQ_1)}{\hat \ttY=\hat \ttP}\,,
\quad
\casi{\ttX_1=\ttQ_1}{\hat \ttX=\hat \ttQ + \hat\qq(\hat\ttP,\ttQ_1)}
 \,,
\eeq
satisfies:
\beq{elficheck}
\Fiduec: \ZZ_{2 \chr}\times \torus_\chs\to \ZZ_{3 \chr}\times \torus_\chs\,,
\eeq
and 
\beqa{hamsecd}
\hamsecd(\ttP,\ttQ_1)&:=&\hamsecu\circ \Fiduec(\ttP,\ttQ_1)\\
&=& \big(1+\tilde\cin(\ttP,\ttQ_1)\big)\, 
\big(\ttP_1-\pp_{\rm o}(\hat \ttP)\big)^2 + \Gfo(\hat\ttP)+ \Gf(\hat\ttP,\ttQ_1)\,,
\nonumber
\eeqa
for suitable functions  $\tilde\cin$, $\Gfo$ and $\Gf$ (explicitly defined in \equ{limone} below, with $\langle\Gf\rangle=0$) real analytic on, respectively,  
 $\ZZ_{2\chr}\times \torus_\chs$, $\hat \ZZ_{2\chr}$  and  $\hat \ZZ_{2\chr}\times \torus_\chs$ , 
which satisfy the bounds:
\beq{cima}\ts
|\tilde \cin|_{2\chr,\chs}\leq \frac{\bettao}{4\chr^2}\,,
\qquad
|\Gfo|_{2\chr}\leq 3 \bettao\,\qquad
\quad
|\Gf-\bGf|_{2\chr,\chs}\leq 2\betta\,.
\eeq
\end{lemma}

\proof 
We start by proving \equ{tettapic}. Recalling   \equ{fangorn}, \equ{bollettino1} and \equ{bollettino3}, we have  
\beq{chenoia}\ts
\frac{\pi+\chs}\chs\stackrel{}\le 1+ 2\pi \ttcs < 8 \ttcs<\K\,.
 \eeq
Now, by the definitions in \equ{betta},   \equ{spigolak}, \equ{cerbiatta}, \equ{dublino}, we find
\beqno\ts
\frac{\bettao}{\chr^2}
=
4  \itcd^2 \, \frac{|k|^2}{\K^{14n+4}}
\stackrel{\equ{chenoia}}\le 4  \itcd^2 \, \frac{1}{\K^{14n+1}} \frac\chs{\pi+\chs}\,,
\eeqno
which yields \equ{tettapic} since, by assumption,  $\K>\KO\ge \itcd$.

\nl
(i) Let us denote by ${\bf X}:= \hat\ZZ_{4 \chr,\chs}\times  \torus_\chs$ and by
 $\mathcal X$  the complete metric space formed by the real analytic complex--valued functions $u: {\bf X}\to \{z\in\complex:|z|\le \chr/2\}$, equipped with the metric given by   the distance in sup--norm on $\bf X$. Let us also denote:
\beq{pontediferro}
\ttG^\sharp:=\Guo+\Gu\,,\qquad\quad \tilde \ttG^\sharp := \Gu- \bGf\,.
\eeq
Note that $\Gu$ and $\tilde \ttG^\sharp$ have zero average.
Consider the operator $F:u\in {\mathcal X} \mapsto F(u)$, where
$F(u)(\hat \ttP,\ttQ_1):=-{\textstyle \frac12} \partial_{\ttY_1} \ttG^\sharp (u(\hat\ttP,\ttQ_1),\hat\ttP, \ttQ_1)$.
If $u\in  {\mathcal X}$, then, by Cauchy estimate we get
\beqa{onam1}
\sup_{\bf X} |F(u)|&=&{\textstyle \frac12} \sup_{\bf X}
\big|  \partial_{\ttY_1} \ttG^\sharp(u(\hat\ttP,\ttQ_1),\hat\ttP, \ttQ_1)\big|
\nonumber\\
&=&{\textstyle \frac12} \sup_{\bf X}
\big|  \partial_{\ttY_1} \big[\ttG^\sharp(u(\hat\ttP,\ttQ_1),\hat\ttP, \ttQ_1)- \bGf(\ttQ_1)\big]\big|
\nonumber\\
&\le &
 \frac12\ \frac{\big|  \ttG^\sharp-\bGf\big|_{4 \chr,\chs}}{4 \chr- \frac\chr2}
\nonumber\\
&\stackrel{\equ{betta}}\le&  \frac12\, \frac{\bettao+\betta }{4 \chr- \frac\chr2} 
\le \frac27\, \frac\bettao{\chr}
\stackrel{\equ{tettapic}}{< } \frac27\, \chr<\frac\chr2\,.
 \eeqa
Thus,  $F:{\mathcal X}\to {\mathcal X}$. 
Let us check that $F$ is, in fact,  a contraction  on  ${\mathcal X}$. If $u,v\in{\mathcal X}$, then, again, by Cauchy estimate, 
\equ{betta} and \equ{tettapic}, we get\footnote{$u$ and $v$, in the r.h.s. of the first inequality, are evaluated at $(\ttQ_1,\hat\ttP)$. }
\beqa{onam2}
\sup_{\bf X} |F(u)-F(v)| &\le & 
{\textstyle \frac12} \sup_{\bf X}
\big|  \partial_{\ttY_1}\big( \ttG^\sharp(u,\hat \ttP,\ttQ_1)-  \ttG^\sharp(v,\hat \ttP,\ttQ_1)\big)\big|
\nonumber\\
&\le&  \frac12\  \big|  \partial_{\ttY_1}^2(\ttG^\sharp-\bGf\big)|_{\frac\chr2,\chs}\cdot \sup_{\bf X} |u-v|
\nonumber\\
&\le &
 \frac12\ \frac{\big|  \ttG^\sharp-\bGf\big|_{4 \chr,\chs}}{\big({4 \chr}- \frac\chr2\big)^2}\cdot \sup_{\bf X} |u-v|
\nonumber\\
 &\stackrel{\equ{betta}}{\le}& \frac{4}{49}\ \frac{\bettao}{\chr^2} \cdot \sup_{\bf X} |u-v|
 \stackrel{\equ{tettapic}}{<} \frac18\ \cdot \sup_{\bf X} |u-v|\,,
 \eeqa
showing that $F$ is a contraction on $\mathcal X$. Thus, by the standard Contraction Lemma, it follows that 
there exists 
a unique $\pp\in {\mathcal X}$ solving \equ{fpe}. 

\nl
Since $\pp=F(\pp)$, one sees that \equ{pitale} follows from \equ{onam1}.\\
The first bound in \equ{tess} follows immediately from \equ{pitale}.
\\
To prove the second estimate in \equ{tess}, write\footnote{To simplify notation, we drop, here, from the notation the explicit dependence on $\hat\ttP$ and $\ttQ_1$ of $\ttG^\sharp$.}
\beq{deangelis}
 \partial_{\ttY_1} \ttG^\sharp (\pp)= 
 \partial_{\ttY_1} \ttG^\sharp (\pp_{\rm o}+\tilde\pp) =
 \partial_{\ttY_1} \ttG^\sharp (\pp_{\rm o}) + w \tilde\pp\,,\quad {\rm with} \quad w:= \int_0^1   \partial^2_{\ttY_1} \ttG^\sharp (\pp_{\rm o}+t \tilde\pp)dt\,.
\eeq
As above, by Cauchy estimates, 
\beq{marcore}
|w|_{4\chr,\chs}\le  \frac2{49} \, \frac{\bettao+\betta}{\chr^2}<\frac18\,.
\eeq
Thus, by \equ{deangelis}, Cauchy estimates, and  \equ{marcore},  observing 
that\footnote{Recall that
$\langle \Gu(\ttY,\cdot)\rangle=0$ as stated in Lemma \ref{phi1}.} $\langle  \partial_{\ttY_1} \Gu(\pp_{\rm o})\rangle=0$, one finds
\beqano
|\tilde \pp|=|\pp-\pp_{\rm o}|&\stackrel{\equ{deangelis}}=&{\textstyle \frac12} \big|  \partial_{\ttY_1} \ttG^\sharp (\pp_{\rm o}) - \langle  \partial_{\ttY_1} \ttG^\sharp (\pp_{\rm o})\rangle +
w \tilde \pp - \langle w \tilde \pp \rangle\big|\\
&=& 
{\textstyle \frac12} \big|  \partial_{\ttY_1} \Gu (\pp_{\rm o}) - \langle  \partial_{\ttY_1} \Gu(\pp_{\rm o})\rangle +
w \tilde \pp - \langle w \tilde \pp \rangle\big|\\
&=& 
{\textstyle \frac12} \big|  \partial_{\ttY_1} \big(\Gu (\pp_{\rm o})-\bGf\big)+
w \tilde \pp - \langle w \tilde \pp \rangle\big|\\
&\stackrel{\equ{betta},\equ{marcore}}\le& \frac12 \Big(\, \frac27 \, \frac\betta\chr\Big) + \frac12 |\tilde \pp|\,,
\eeqano
which yields immediately the second bound in \equ{tess}.
\\
Next, since 
$\tilde \pp$ has zero average over the torus, the function $\phi$ defined in \equ{guaito}
defines a (real analytic) periodic function such that $\partial_{\ttX_1}\phi=\tilde \pp$.
Furthermore, by the second estimate in  \equ{tess}, one has\footnote{$\pi+\chs$ is an estimate of the length of the integration path in \equ{guaito}, as the real part of
 $\ttQ_1$ can be taken in $[-\pi,\pi)$.}
$
|\phi|_{4 \chr,\chs}< \frac{\betta}{3 \chr}\ (\pi+\chs)\,,
$
so that, by Cauchy estimates, also  last bounds in \equ{tess} follow.

\nl
(ii)  By the definition of $\Fidue$ in \equ{Fidue},  by \equ{tess} and \equ{tettapic}, the relations in
\equ{elficheck} follow at once.
\\
Now, define\footnote{Here, $\pp_{\rm o}=\pp_{\rm o}(\hat\ttP)$.}
\beq{limone}
\left\{\begin{array}{l}
\dst \tilde\cin(\ttP,\ttQ_1):=\int_0^1 (1-t)\partial_{\ttY_1}^2 \ttG^\sharp\big(\pp_{\rm o}+t (\ttP_1-\pp_{\rm o}),\hat\ttP,\ttQ_1\big)dt\,,\\
\Gfo(\hat\ttP):= \langle \pp(\hat \ttP,\cdot)^2\rangle +\langle \ttG^\sharp\big( \pp(\cdot,\hat \ttP),\hat\ttP,\cdot)\rangle
\,,\\
\dst \Gf(\hat\ttP, \ttQ_1):= \pp(\hat \ttP,\ttQ_1)^2
+ \ttG^\sharp\big( \pp(\ttQ_1,\hat \ttP),\hat\ttP,\ttQ_1) - \Gfo(\hat\ttP)
\,,
\end{array}\right.
\eeq
then, 
by Taylor's formula, \equ{hamsecu}, \equ{Fidue}, \equ{pontediferro}  and \equ{fpe}, one finds\footnote{$\dst g(t_0+\t)=g(t_0)+g'(t_0)\t+\Big(\int_0^1 (1-t) g''\big(t_0+t\t\big)dt\Big) \t^2$ with $g=\Gu(\cdot,\ttQ_1)$, $t_0=\pp$ and $\t=\ttP_1-\pp_{\rm o}$. For ease of notation we drop the (dumb) dependence upon $\hat \ttP$ in these formulae.}
\beqa{onam3}
\hamsecd(\ttP_1,\ttQ_1)&:=&\hamsecu\circ \Fiduec(\ttP_1,\ttQ_1)
= (\ttP_1+\tilde\pp)^2+\ttG^\sharp(\ttP_1+\tilde\pp,\ttQ_1)\nonumber\\
&\stackrel{\equ{guaito}}=&\big(\pp+(\ttP_1-\pp_{\rm o})\big)^2+ \ttG^\sharp\big(\pp+(\ttP_1-\pp_{\rm o}),\ttQ_1\big)\nonumber\\
&\stackrel{\equ{limone}}{=}&  (\ttP_1-\pp_{\rm o})^2+2(\ttP_1-\pp_{\rm o})\pp + \pp^2+ \ttG^\sharp(\pp,\ttQ_1)+\partial_{\ttY_1}\ttG^\sharp(\pp,\ttQ_1) (\ttP_1-\pp_{\rm o}) \nonumber \\
&& + (\ttP_1-\pp_{\rm o})^2 \tilde \cin\nonumber\\
&\stackrel{\equ{fpe}}{=}&
(1+\tilde \cin) (\ttP_1-\pp_{\rm o})^2 +  \pp^2+ \ttG^\sharp(\pp,\ttQ_1)\nonumber\\
&\stackrel{\equ{limone}}{=}& (1+\tilde \cin) (\ttP_1-\pp_{\rm o})^2+ \Gfo+  \Gf(\ttQ_1)\,,
\eeqa
proving \equ{hamsecd}.
\\
Let us now prove \equ{cima}.
Observe that for $\ttP\in \ZZ_{2\chr}$  by \equ{tettapic} and \equ{tess}
the segment  $\big(\pp_{\rm o}+t (\ttP_1-\pp_{\rm o}),\hat \ttP\big)$, $t\in[0,1]$, 
still belongs to $\ZZ_{2\chr}$, hence, by
definition of $\tilde \cin$ in \equ{limone}, by Cauchy estimate\footnote{Compare, also, the estimates done in \equ{onam2}.} and \equ{betta} one obtains the first
estimate in \equ{cima}.
\\
By the definition of $\Gfo$, by \equ{pitale} and \equ{betta},   observing that $|\ttG^\sharp|\le \bettao+\betta\le 2\bettao$, one gets immediately the second estimate in \equ{cima}.
\\
As for the third estimate in \equ{cima}, by the definitions given, one has that\footnote{Dropping, again,  in the notation the dumb variable $\hat \ttP$.}
\beq{proietti}
\Gf-\bGf
=
\big(\pp^2-\langle \pp^2\rangle\big)+\big(\ttG^\sharp(\pp,\cdot) - \langle \ttG^\sharp(\pp,\cdot)\rangle-\bGf\big) \,.
\eeq
Let us estimate the terms in brackets separately. For $\hat\ttP\in \hat\ZZ_{2\chr}$ and $\ttQ_1\in\torus_{\chs}$, one finds
\beqa{fiorini}
|\pp^2-\langle \pp^2\rangle|=|2\tilde \pp\pp_{\rm o} +\tilde\pp^2-\langle \tilde \pp^2\rangle|\le
(2|\pp_{\rm o}|+2|\tilde \pp|)\, |\tilde\pp| 
\stackrel{\equ{tess}}{\le} 
\frac49 \bettao \frac\betta{\chr^2}\stackrel{\equ{tettapic}}{<} \frac\betta2\,.
\eeqa
To estimate the second term in \equ{proietti}, we define
$$\z(t):=
\ttG^\sharp(\pp_{\rm o}+t \tilde \pp,\ttQ_1) - \langle \ttG^\sharp(\pp_{\rm o}+t \tilde \pp,,\cdot)\rangle\,,
$$
and observe that (recall \equ{pontediferro})
$\z(0)= \Gu(\pp_{\rm o},\ttQ_1)$
and that, by Cauchy estimates, we get\footnote{Reasoning as in \equ{onam1}.}
\beqa{oji}
|\z'(s)|&\le& |\tilde \pp|\,  \int_0^1\big| \partial_{\ttY_1}\big( \ttG^\sharp(\pp_{\rm o}+t \tilde \pp,\ttQ_1) - \langle \ttG^\sharp(\pp_{\rm o}+t \tilde \pp,,\cdot)\rangle\big)\big|dt\nonumber\\
&\le& |\tilde \pp|\,    \frac{2|\ttG^\sharp|_{4\chr ,\chs}}{4\chr -\frac\chr2}
\le |\tilde \pp|\,    \frac{4\bettao}{4\chr -\frac\chr2} 
\stackrel{\equ{tess}}{<}\frac8{21} \frac{\bettao}{\chr^2}\ \, \betta\stackrel{\equ{tettapic}}{<}\frac12 \betta\,.
\eeqa
Thus, 
\beqno
\big|\ttG^\sharp(\pp,\cdot) - \langle \ttG^\sharp(\pp,\cdot)\rangle-\bGf\big|\le |\z(0)-\bGf|+\int_0^1|\z'(t)|dt
\le
| \Gu(\pp_{\rm o},\cdot)-\bGf|+\frac12 \betta\stackrel{\equ{betta}}{\le} \frac32 \betta\,.
\eeqno
Putting together this estimate and \equ{fiorini} one gets also the third estimate in \equ{cima}. \qed

\nl
The final   transformation is again just a translation, which is done so that
 {\sl all equilibria of the
secular system will lie on the angle--axis in its  2--dimensional phase space}. 

\begin{lemma}\label{platano}
The  real analytic symplectic map $\Fitre\in\Gdag$ defined as
\beq{Fitre}
\Fitre:
(\ttp,\ttq) \mapsto
(\ttP,\ttQ)
 \quad{\rm with}\quad
\casi{\ttP_1=\ttp_1 + \pp_{\rm o}(\hat\ttp)}{\hat \ttP=\hat \ttp}\,,
\qquad
\casi{\ttQ_1=\ttq_1}{\hat\ttQ=\hat\ttq -\ttq_1\partial_{\hat \ttp} \pp_{\rm o}(\hat\ttp) \, }
 \,,
\eeq
satisfies
\beq{giovannino}\textstyle
\Fitrec: \ZZ_{\chr}\times \torus_\chs\to \ZZ_{2\chr}\times \torus_\chs\,.
\eeq
Furthermore, one has:
\beq{Hsharp-}
\hamsecd\circ \Fitrec(\ttp,\ttq_1)= 
\big(1+\cins(\ttp,\ttq_1)\big)\,  \ttp_1^2 + \Gfo(\hat\ttp)+  \Gf(\hat\ttp,\ttq_1)\,,
\eeq 
where 
\beq{cins}
\cins(\ttp,\ttq_1):= \tilde\cin\big(\pp_{\rm o}(\hat\ttp)+\ttp_1, \hat\ttp,\ttq_1\big)\,,
\eeq
and the following bunds hold:
\beq{cimadue}
|\cins|_{\chr,\chs}\leq \frac{\bettao}{4\chr^2}\,,
\qquad
|\Gfo|_{2\chr}\leq  3 \bettao\,\qquad
\quad
|\Gf-\bGf|_{2\chr,\chs}\leq  2\betta\,.
\eeq
\end{lemma}
\proof  Just observe that, if $|\ttp_1|<\chr$, then, by \equ{tess} and \equ{tettapic}, it follows that,
for all $\ttp\in\ZZ_{\chr}$,
$$\ts
|\pp_{\rm o}(\hat\ttp)+\ttp_1|<  \frac{\bettao}{3 \chr}+\chr\le \frac\chr3+\chr=\frac43 \chr<2 \chr\,,
$$
so that \equ{giovannino} holds.
Finally, by \equ{cima}, we get\footnote{{\sl $\Gf$ and $\bGf$ are the same} as in \equ{cima} of Lemma~\ref{avogado}.} \equ{cimadue}. \qed
We are ready for the 

\nl
\proof {\bf of Theorem \ref{sivori}}\\
Recall the definitions of $\Phi_{\!{}_j}$, $1\le j\le 3$,  in, respectively,  Lemma~\ref{phi1}, \equ{Fidue} and \equ{Fitre} and define
$\Phi_\diamond:=\Fiuno\circ\Fidue\circ \Fitre$, 
and
$\hzk(\hat\p):=
\frac{1}{|k|^2} 
\normadue \proiezione^\perp_k \hAA^T \hat\p\normadue ^2
+\Gfo(\hat\p)$. Then, the expression for $\Hsharp$ in \equ{tikitaka}
 follows by
\equ{spigolak}, \equ{hamsecd} and \equ{Hsharp-}.
\\
By \equ{giovannino} and  Lemma~\ref{platano},
the Hamiltonian function $\Hsharp$ is real analytic on $\ZZ_\chr\times\torus_\chs$, where  
$\ZZ=(-\Ro,\Ro)\times\hat \ZZ$
(compare \equ{cerbiatta}).
\\ 
By \equ{P2+} and Proposition~\ref{punti} we have 
that $\bGf$ in \eqref{paranoia}  is $\morse$--Morse with $\morse$ as in \equ{cerbiatta}.
\\
Let us, now, estimate $|\bGf|_\so$. Consider, first, $\noruno{k}<\Nf$. Then, estimating $|f_{jk}|$ by $e^{-|j|\noruno{k}s}$,  by the definition of $\so$, we get
\beqano
|\bGf|_\so&\stackrel{\equ{paranoia}}=& \e_k |\fproj f|_{\so}
\le 
\e_k |\fproj f|_{s/2}=\e_k \sum_{j\neq 0} |f_{jk}|e^{\frac{|j|\noruno{k}s}{2}} 
\\
&\le& \frac{8\e}{|k|^2} \, \frac{e^{-s/2}}{2(1-e^{-s/2})}<  \frac{8\e}{|k|^2} \, \frac{1}{s}\,.
\eeqano
If $\noruno{k}\ge\Nf$ one has
\beqno
|\bGf|_\so=|\bGf|_1
\stackrel{\equ{alfacentauri}}=\frac{4\e}{|k|^2}|f_k||\cos(\sa+ \sa_k)+F^k_\star(\sa)|_1
\stackrel{\equ{gallina}}{\le}\frac{4\e}{|k|^2} \, |f_k|\, (\cosh 1+2^{-40})<\frac{8\e}{|k|^2} |f_k|\,.
\eeqno
Thus, by definitions of $\chi_{{}_k}$ in \equ{cerbiatta} and $\bttcd$, one gets 
\beq{crusca}
|\bGf|_\so\le \suca \,,
\eeq
with $\suca$ as in \equ{cerbiatta}.
Next, since  $\chi_{{}_k}\le 1\le \bttcd$, 
\beq{mappo}
|\Gf-\bGf|_{\ro,\so}\stackrel{\equ{cima}}\le 2\betta\stackrel{\equ{betta}}=\frac{8\e}{|k|^2} \chi_{{}_k}\tetta
\stackrel{\equ{cerbiatta}}\le \suca \lalla\,.
\eeq
By \equ{cimadue}, \equ{cerbiatta}, \equ{betta}, using the inequalities $|k|\le \KO\le \K/6$,
recalling \equ{tikitaka}, \equ{dublino}, and the hypothesis $\KO\ge \itcd$ (in the last inequality), 
 one sees that 
\beq{mappo2}
|\cins|_{\ro,\so}\leq  \itcd^2\, \frac{|k|^2}{\K^{2\nu}}\, \tetta\le 
\frac{\itcd^2}{36}\, \frac{1}{\K^{2(\nu-1)}} \, \tetta
<
\tetta=\lalla\,.
\eeq
Then \equ{cimabue}
follows by \equ{crusca}, \equ{mappo} and
\equ{mappo2}.
\\
Finally, observe that, by the definitions in
\equ{cerbiatta} and \equ{betta}
one has
\beq{stent}\ts
{\suca}/{\morse}=
\casitwo{ \frac{4 \bttcd}{\b}
\,,}{\noruno{k}<\Nf\,,}{4 \bttcd
\,,\phantom{\dst \int}
}{\noruno{k}\ge \Nf\,.}
\eeq
Then, \equ{alce} with $\varpi$ as in \equ{kappa}, follows immediately by the definitions in 
\equ{cerbiatta}, 
\equ{sucamorse} and \equ{stent}.  
\qed

\small


%
%
%
%
%

\end{document}